\documentclass[review]{elsarticle}

\usepackage{lineno,hyperref}
\usepackage{float, comment}

\usepackage{amsmath,amsxtra,latexsym,amsthm, amssymb, amscd}
  \usepackage{epsfig} 

\usepackage{graphicx}
  \usepackage{epstopdf}

\modulolinenumbers[5]

\newtheorem{remark}{Remark}

\begin{document}

\begin{frontmatter}

\title{A stochastic differential equation model for foraging behavior of fish schools}


\author[Ton]{T$\hat{\text o}$n Vi$\hat{\d{e}}$t T\d{a}}        
\ead{tavietton@agr.kyushu-u.ac.jp}

\author[Linh]{Linh Thi Hoai Nguyen}
\ead{nth.linh@ist.osaka-u.ac.jp}

\address[Ton]{Center for Promotion of International Education and Research\\
Faculty of Agriculture, Kyushu University\\
6-10-1 Hakozaki, Higashi-ku, Fukuoka 812-8581, Japan}

\address[Linh]{Department of Information and Physical Sciences\\
Graduate School of Information Science and Technology, Osaka University\\
1-5 Yamadaoka, Suita, Osaka 565-0871, Japan}

\begin{abstract}
We present a novel model of stochastic differential equations for foraging behavior of fish schools  in space including   obstacles. We then study the model numerically. Three configurations of space with different locations of  food resource are considered. In the first configuration, fish move in free but limited space. All individuals can find food almost surely. In the second and third configurations, fish move in limited space with one or two obstacles. 
Our results reveal that on one hand, when school size  increases, so does the probability of foraging success.  On the other hand, when it exceeds an  optimal value, the probability   decreases.  In all configurations, fish always keep a school structure through the process of foraging.
\end{abstract}



\end{frontmatter}

\allowdisplaybreaks


\section{Introduction}  \label{introduction}
Swarming of animals, for example schooling of fish, flocking of birds, or herding of mammals, is one of the most commonly observed phenomenon in the real world but a challenge  to study in biology. Swarm behavior is a collective behavior exhibited by animals of similar size which aggregate together. This remarkable phenomenon has already attracted interest of researchers from diverse fields including biology, physics, mathematics, and computer engineering.

In order to understand swarming dynamics, one way is to construct mathematical models on the basis of local rules, which are observed by researchers in discipline. 
Vicsek et al. (\cite{Vicsek1995}) modeled the movement of self-driven particles by difference equations in which  particles move at constant speed and choose their new heading to be the average of those of nearby particles located within a unit distance. Based on this model, Cucker-Smale (\cite{Cucker2}) presented a simple model of  ordinary differential equations  for swarming by using an interaction between individuals  (i.e., alignment). Some stochastic versions of the Cucker-Smale model are studied in \cite{Cucker1,Ton}.

Oboshi et al. (\cite{Oboshi2002}) modeled  schooling by difference equations setting a rule that each individual choses one way of action among four possibilities according to a distance to the closest mate. Olfati-Saber (\cite{RezaMultiAgent}) and D'Orsogna et al. (\cite{DOrsogna2006}) presented differential equation models, but deterministic ones, utilizing the generalized Morse function and attractive/repulsive potential functions, respectively. Gunji et al.  (\cite{Gunji1999}) considered dual interaction which produced territorial and schooling behavior. 

In the above researches, swarming is considered in free spaces. In \cite{Gautrais2}, Gautrais et al.  characterized the spontaneous behavior of a single fish ({\it Kuhlia mugil}) in limited spaces (shallow circular swimming pools). In this model, the fish moves at constant speed but the angular velocity of the fish orientation obeys  a stochastic differential equation (or a Ornstein-Uhlenbeck process) which asserts that  the fish can avoid collisions with the tank walls. Gautrais et al. (\cite{Gautrais1}) then used a bottom-up methodology to construct models of group motion from data gathered at an individual scale.
By analyzing experimental  data captured from individual zebrafish  ({\it Danio rerio}) via automated visual tracking, Zienkiewicz et al. (\cite{Zienkiewicz}) suggested that the model proposed by Gautrais et al. in which speed is constant may not suitable to describe the single and collective locomotion of zebrafish. The authors then extended the approach by Gautrais et al. to construct a data-driven model for locomotion of a single zebrafish. In this model, both the speed and the angular velocity of the fish are Ornstein-Uhlenbeck processes.  For further discussion about how and why animals interact, see  \cite{Danchin,Sumpter1,Sumpter2,Vicsek2012} and references therein.

In  the  monograph \cite{Camazine2001}, Camazine et al. presented an insight of picking up the following behavioral rules of individual fish:
\begin{enumerate}
\item [(a)] The school has no leaders and each fish follows the same behavioral rules.
\item [(b)] To decide where to move, each fish uses some form of weighted average of the position and orientation of its nearest neighbors.
\item [(c)] There is a degree of uncertainty in the individual's behavior that reflects both the imperfect information-gathering ability of a fish and the imperfect execution of the fish's actions.
\end{enumerate}

We should mention that these three rules are  based on empirical results of Aoki (\cite{Aoki1982}), Huth-Wissel (\cite{Huth1992}) and Warburton-Lazarus (\cite{WL}). And similar assumptions, but deterministic ones, were also introduced by Reynolds (\cite{Reynolds1987}). 

Based on this idea, the authors of the present paper published a few papers in the view point of mathematical science. In Uchitane-T\d{a}-Yagi (\cite{Uchitane2012}), we used stochastic differential equations (SDEs) to construct  a mathematical model  describing  the process of schooling of $N$-fish system in a non-limited space, i.e., the Euclidean space $\mathbb R^d \,(d=1,2,3\dots).$

In Nguyen-T\d{a}-Yagi (\cite{Nguyen2014}), we gave quantitative investigations for that model. 
 We  performed  numerical computations to clarify some important effects of parameters of the model on determining geometrical structures of school. 

In Nguyen-T\d{a}-Yagi (\cite{LinhTonYagi}), we constructed a mathematical model of SDEs which describe the movement of individuals in space with obstacle. To construct that model, we introduced a local rule of obstacle avoidance for individual fish:
\begin{enumerate}
\item[(d)] Each fish executes an action for avoiding obstacle according to the reflection law of velocity with a weight depending on distance.
\end{enumerate}
By numerical computations, we found four obstacle avoidance patterns of fish school, named  {\it Rebound},  {\it Pullback},  {\it Pass and Reunion}, and   {\it Separation}. In addition, we showed how these patterns change as  crucial modeling parameters  change.

We are now interested in foraging behavior of fish schools  in noisy environment. Previous experimental observations showed that  swarming is beneficial to foraging. Gotmark et al. (\cite{Gotmark}) showed that the foraging success of gulls ({\it Larus ridibundus})  increases with flock size up to at least eight birds. Couzin et al. (\cite{Couzin1})  performed  experiments  in a shallow tank on school of $2^n$- golden shiners fish ({\it Notemigonus crysoleucas}) $(n=1,2\dots 6)$ in which fish track the preferred, darker regions  of a circular patch  (darkest at its center and transitioned to the brightest light levels) that move at a constant speed  in the tank. It is shown that when school size increases, so does school-level responsiveness to the environment. In other words,  large schools track target better than smaller schools.

In order to simulate collective foraging, Shklarsh et al. (\cite{Shklarsh}) used local rules of individuals in swarm  (repulsion, attraction, alignment, and reaction to the environment) (\cite{Couzin2,Couzin3}) to  construct a model  of {\it difference equations} for collective navigation of bacteria-inspired smart agents in complex terrains. The authors showed that the  length of path (from starting point to a fixed target) decreases as a function of group size  due to collection of information from more agents. In the movement of agents to the target, the group may {\it separate into many clusters}.

In the present paper,  we  present a model of SDEs for foraging behavior of fish schools in noisy environment with  obstacles. For this purpose, we newly introduce a  local rule of individuals for foraging. We then write out the rule  into a mathematical formula, and integrate it into our previous models to obtain the desired model. As a consequence, individuals always keep school structure (i.e., the group is not separated into clusters) through the process of foraging. We then numerically study the model in three configurations of obstacle with different locations of food resource. Our numerical results qualitatively agree with the above experimental evidences  (\cite{Couzin1,Gotmark}) or empirical results  (\cite{Pitcher})  that the bigger the school size the larger the probability of foraging success. 

The increase in probability of foraging success is, however, not retained unboundedly.  As school size exceeds some optimal value, the probability decreases. Our model therefore may give an estimate for that optimal size for each species.

The organization of  the paper is as follows. In Section \ref{sec2}, we introduce our SDE models. Subsections \ref{sec2.1} and \ref{sec2.2}  review our previous SDE models in free spaces and in spaces with obstacle. Subsection \ref{sec2.3} presents  the local rule for foraging and integrate it into the previous models. The resulting model of SDEs then describes foraging behavior of fish schools in noisy environment with obstacles. 
 Section \ref{sec3} gives numerical results which agree with empirical and experimental evidences.  The paper ends with some conclusions in Section \ref{section5}.

\section{Mathematical models for fish schooling}\label{sec2}

\subsection{SDE model in free spaces} \label{sec2.1}

Recently, we introduced a SDE model for fish schooling based on the three behavioral rules {\rm (a)-(c)}  in the Introduction (\cite{Uchitane2012}). Each of $N$ fish is regarded as a  particle moving in the free space $\mathbb R^d$ $(d=1,2,3\dots)$. 
The interactions between particles in our model include generalization of the inverse-square law of universal gravitation (attraction), generalization of the Van der Waals forces (repulsion), and the alignment of particles. 
 Our SDE model  reads as
\begin{equation}\label{eq0}
\begin{aligned}
\begin{cases}
d{x}_i(t)=&{v}_idt+\sigma_idw_i(t),\\
d{v}_i(t)=&\Big\{-\alpha\sum\limits_{j=1,j\ne i}^N\left(\frac{r^p}{\|{x}_i-{x}_j\|^p}-\frac{r^q}{\|{x}_i-{x}_j\|^q} \right)({x}_i-{x}_j)\\
  &\quad -\beta\sum\limits_{j=1,j\ne i}^N\left(\frac{r^p}{\|{x}_i-{x}_j\|^p}+\frac{r^q}{\|{x}_i-{x}_j\|^q} \right)({v}_i-{v}_j)+F_i({x}_i,{v}_i)\Big\}dt.
\end{cases}
\end{aligned}
\end{equation}
Here, ${x}_i(t)$ and ${v}_i(t)$ $(i=1,2\dots N)$ denote the position and the velocity, respectively, of the $i$-th individual at time $t$; and $\|\cdot\|$ denotes the Euclidean norm of vector.

 The first equation of \eqref{eq0} is a stochastic equation for the unknown ${x}_i(t)$, where $\sigma_idw_i$  denotes a stochastic differentiation of a $d$-dimensional independent Brownian motion defined in a filtered probability space.  The second one is a deterministic equation for the unknown ${v}_i(t)$, where $1<p<q<\infty$ are fixed exponents; $\alpha$ and $\beta$ are positive coefficients of attraction and of velocity matching among individuals, respectively;  $r>0$ is a fixed number; $F_i({x}_i,{v}_i)$ stands for an external force function acting on the $i$-th individual.

If the $i$-fish is far from the $j$-fish, i.e., $\|{x}_i-{x}_j\|>r,$ then it would move towards the $j$-th due to the attraction force. To the contrary, if they are close enough, i.e., $\|{x}_i-{x}_j\|<r,$ then they would avoid collision with each other due to the repulsive force. The quantity $r$ therefore plays as the critical distance. 

The velocity matching of the $i$-th individual to the $j$th individual also has a similar weight depending on the distance $\|{x}_i-{x}_j\|$. Degree of matching is higher when $\|{x}_i-{x}_j\|<r$ than the contrary case in order to avoid collisions.

In the meantime, the exponent $p$ denotes a degree of how far the attraction reaches. If $p$ is large, then the attraction range is small. If $p$ is small, then individuals can attract each other even if there is a long distance among them. 

An advantage of using SDE models like \eqref{eq0} may be the easiness of mathematical treatments. One can utilize the well-developed theory of SDEs and the numerical methods (\cite{Kloeden2005}). Its flexibility may be another advantage. As seen in the next two subsections, we can make new models  simply by introducing suitable external force functions $F_i$ in \eqref{eq0}.

\subsection{SDE model in spaces with obstacle}  \label{sec2.2}

Our SDE model in spaces with obstacle was presented in  \cite{LinhTonYagi}, in which the obstacle is a compact sphere. We gave a specific form for the external force function $F_i$ in the system  \eqref{eq0} in order to describe obstacle avoidance of fish on the basis of the local rule {\rm (d)} (see Section \ref{introduction}). Let us review that model with a slight change of the shape of obstacle. 

Denote by $S$ the surface of a static obstacle, say a long thin rectangular parallelepiped. Assume that the $i$-th particle is at position $x_i,$ and has velocity $v_i$. 
The particle avoids $S$ by matching its velocity to the reflection vector of $v_i,$ say $\text{\rm Rf}(x_i,v_i; S)$ with respect to $S$. This reflection vector is defined as follows.

\begin{figure}[ht]
 \begin{center}
\includegraphics[width=13cm, height=6cm]{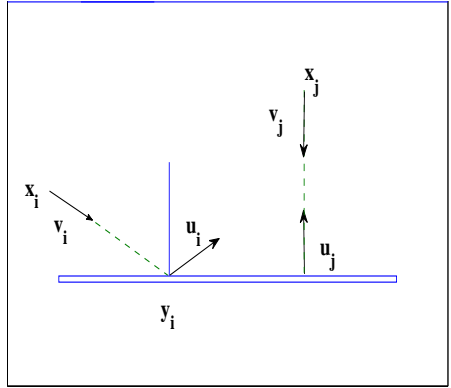}
 \caption{Rule of obstacle avoidance in 2-dimensional space. A fish at position $x_i$ (resp. $x_j$) and spees $v_i$ (resp. $v_j$) avoids collision with obstacle by matching its velocity to its reflection vector $u_i$ (resp. $u_j$).}
 \label{Fig1}
 \end{center}
 \end{figure}

 Let $l_i$ be the ray with origin $x_i$ and direction $v_i$, i.e., 
$$l_i=\{x\in \mathbb R^d\,{;}\,x=x_i+s v_i, 0\leqslant s <\infty\}.$$
If $l_i$ meets $S$ at a point $y_i\in S$, we define $\text{\rm Rf}(x_i,v_i;S)$ to be the reverse of the symmetrical vector of $v_i$ to the line which is perpendicular to $S$ at $y_i$. In the particular case where   $l_i$ is perpendicular to $S$, we observe that  $\text{\rm Rf}(x_i,v_i;S)=-v_i$. If $l_i$ does not meet $S$ (including the case of $v_i=0$), we define $\text{\rm Rf}(x_i,v_i;S)=v_i$. Figure \ref{Fig1} illustrates $(x_i,v_i)$ and $u_i=\text{\rm Rf}(x_i,v_i;S)$ in two-dimensional space. 

\vspace{0.3cm}

Let  $\Omega \subset \mathbb R^d$ be a domain  for possible movement of all fish. Assume that the boundary $\partial \Omega $ of $\Omega$  contains walls. As above,     $\text{\rm Rf}(x_i,v_i;\partial\Omega)$ then denote the reflection vector of vector $v_i$ starting from $x_i$ with respect to $\partial \Omega$.  
Therefore, each particle is affected locally from the boundary $\partial \Omega$ through a small pinpoint part of boundary at each moment.  

 Analogously to the velocity matching, local affection of $\partial \Omega$ to the motion of the $i$-th fish is then described by the force
 $$F_i({x}_i,{v}_i)= -\gamma \left( \frac{R^P}{\|x_i-y_i\|^P} +
       \frac{R^Q}{\|x_i-y_i\|^Q} \right)[v_i-\text{\rm Rf}(x_i,v_i;\partial \Omega)].$$
Here, $y_i\in\partial \Omega$ is a point at which the fish is feared to collide; $R>0$ is a fixed distance; $\gamma>0$ is a constant; and  $1<P<Q<\infty$ are exponents. 

By this force, if the fish is far from the boundary, i.e., $\|x_i-y_i\|>R,$ then its reaction to avoid obstacle is weak. To the contrary, if the fish is close to the boundary, i.e., $\|x_i-y_i\|<R,$ then it would react more promptly in order to match its velocity to $\text{\rm Rf}(x_i,v_i;\partial \Omega)$. If the ray $l_i$ and $\partial \Omega$ do not intersect, the fish would not take any reaction to the obstacle.

By the above, the SDE model in space including obstacles has the form:
\begin{equation} \label{eq1}
\begin{cases}
  dx_i(t) = & v_i dt+\sigma_idw_i(t),    \\
  dv_i(t) = &\Bigl[  - \alpha \sum\limits_{j=1,\, j\ne i}^N
    \left(\dfrac{r^p}{\|x_i - x_j\|^p}- \dfrac{r^q}{\|x_i-x_j\|^q} \right)(x_i-x_j)   \\
  &\quad - \beta \sum\limits_{j=1,\, j\not=i}^N
    \left(\dfrac{r^p}{\|x_i-x_j\|^p} +\dfrac{r^q}{\|x_i-x_j\|^q} \right)(v_i-v_j)\\
     &\quad -\gamma\left(\dfrac{R^P}{\|x_i-y_i\|^P}+\dfrac{R^Q}{\|x_i-y_i\|^Q}\right)[v_i-\text{\rm Rf}(x_i,v_i;\partial\Omega)]\\
     &\quad +G_i(x_i,v_i)\Bigl]dt,
 \end{cases}
\end{equation}
where $G_i$ denotes again an external force function acting on the $i$-th individual. In \cite{LinhTonYagi}, the case where $G_i\equiv 0$ for all $i$ was studied.

\subsection{SDE model in spaces with obstacle and food resource}   \label{sec2.3}

In this subsection, we newly introduce  a  model of SDEs for foraging behavior of fish schools in noisy environment with  obstacle and food resource.  Consider a fish school moving in a free or limited space to forage for food. The position of food resource is fixed in the space. Fish and food may be separated by obstacles in the sense that the school cannot move to the food in a straightforward way. 

To construct the model, we present here a local rule for foraging:
\begin{enumerate}
\item[(e)] Each fish is sensitive to the gradient of potential formed by scent which is emitted by food, and has tendency to move into a higher direction.
\end{enumerate}
This local rule is then integrated into the model  \eqref{eq1} by giving a specific form of the functions $G_i$ in the second equation of \eqref{eq1}.
\subsubsection{Mathematical formulation of the local rule {\rm (e)}}

Let us make a mathematical formulation of  the local rule  {\rm (e)} for  foraging by using a method of potential functions. 

Let $f$ be the density function of food resource defined in a domain $ \Omega \subset \mathbb R^d$.  Consider an elliptic equation
in $\Omega$ under the homogeneous Neumann boundary condition on $\partial \Omega$: 
\begin{equation}\label{eq2}
\begin{aligned}
\begin{cases}
-c\Delta U+a U=f(x), &  \hspace{2cm}  x\in \Omega, \\
\dfrac{\partial U}{\partial \text{\bf n}}=0, & \hspace{2cm} x \in \partial \Omega.
\end{cases}
\end{aligned}
\end{equation}   
 Here, $U(x)$ denotes the density of scent emitted by food at $x\in \Omega$. The operator $\Delta$ is the Laplace operator in $\Omega$; $c>0$ is a diffusion constant;  $a >0$ is a declining rate of $U(x)$; and  $\text{\bf n}$ denotes the (typically exterior) normal to the boundary $\partial \Omega$.   

We regard $U$  as a potential function. Assume that the $i$-th fish is at  position $x_i \in \Omega$ at some moment.  We choose the external force $G_i$ in \eqref{eq1} to be the gradient of the potential function, i.e., 
\begin{equation} \label{eq3}
G_i(x_i,v_i)=k\nabla U(x_i), \hspace{1cm} i=1,2\dots N,
\end{equation}
where $k>0$ is some sensitivity constant. 
($G_i$ does not depend on velocity $v_i$.)  The  boundary condition in \eqref{eq2} ensures that the domain is perfectly insulated, i.e., scent of food cannot pass through the boundary of the domain.

\subsubsection{Model equations}

We are now ready to state model equations to describe foraging behavior of fish schools in noisy environment with  obstacle and food resource.  The equations are as in  \eqref{eq1}, where $G_i \, (i=1,2\dots  N)$  are defined by \eqref{eq2} and \eqref{eq3}. For convenience, we rewrite the equations as 
\begin{equation} \label{eq4}
\begin{cases}
\begin{aligned}
  dx_i(t) = & v_i dt+\sigma_idw_i(t),   \hspace{1cm} i=1,2\dots N, \\
  dv_i(t) = &\Bigl[  - \alpha \sum\limits_{j=1,\, j\ne i}^N
    \left(\dfrac{r^p}{\|x_i - x_j\|^p}- \dfrac{r^q}{\|x_i-x_j\|^q} \right)(x_i-x_j)   \\
&\quad - \beta \sum\limits_{j=1,\, j\not=i}^N
    \left(\dfrac{r^p}{\|x_i-x_j\|^p} +\dfrac{r^q}{\|x_i-x_j\|^q} \right)(v_i-v_j)\\
&  \quad -\gamma\left(\dfrac{R^P}{\|x_i-y_i\|^P}+\dfrac{R^Q}{\|x_i-y_i\|^Q}\right)[v_i-\text{\rm Rf}(x_i,v_i;\partial\Omega)]\\
 & \quad +k\nabla U(x_i)\Bigl]dt,  \hspace{1cm} i=1,2\dots N,
\end{aligned}\\
\begin{aligned}
&-c\Delta U+a U= f(x),    &  \hspace{2cm}x\in \Omega, \\
&\dfrac{\partial U}{\partial \text{\bf n}}= 0,  &   \hspace{2cm}x \in \partial \Omega.
\end{aligned}
 \end{cases}
\end{equation}

\section{Numerical results}\label{sec3}

In this section, we numerically calculate  the probability of foraging success in the model \eqref{eq4} with different obstacle configurations and locations of food resources, and give some patterns of foraging. Simultaneously, we check the following hypotheses:
\begin{enumerate}
  \item  Fish keep school structure while foraging.
  \item  The experimental observations  (\cite{Couzin1,Gotmark}) and empirical results  (\cite{Pitcher})  say that the bigger the school size the larger the probability of foraging success. 
  \item Optimal value of school size  exists in the sense that the probability of foraging success is highest at that size.
\end{enumerate}

We consider three configurations of space for the model \eqref{eq4}. The first one is a collective foraging scenario of fish schools in a free but limited space (no obstacles). The second and third configurations are the same but we put one or two obstacles inside the space. Positions of food resources in configurations are different.

For simplicity, we  only perform simulations in two-dimensional space.

\subsection{Configuration I}  

Let us consider the system  \eqref{eq4} in a free but limited space, say the rectangle domain $\Omega=[0, 7]\times [0, 4]\subset \mathbb R^2$. Although there isn't any obstacles inside $\Omega$, the term including reflection vector in the right-hand side of the second equation of \eqref{eq4} is still valid. The reason is that if a fish swims towards the boundary of $\Omega$, it would avoid a collision by matching its velocity to its reflection vector.

\begin{figure}[H]
\begin{center}
\includegraphics[width=6cm, height=5cm]{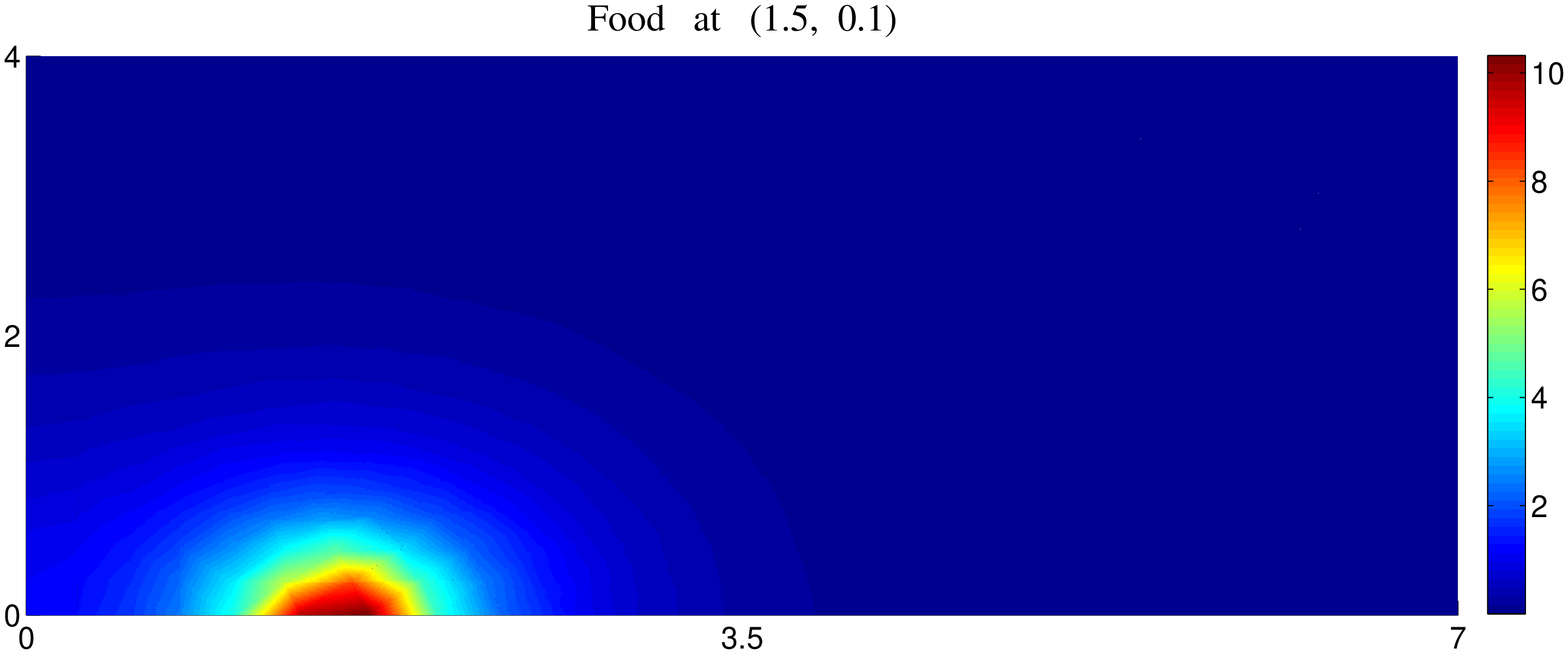} 
\includegraphics[width=6cm, height=5cm]{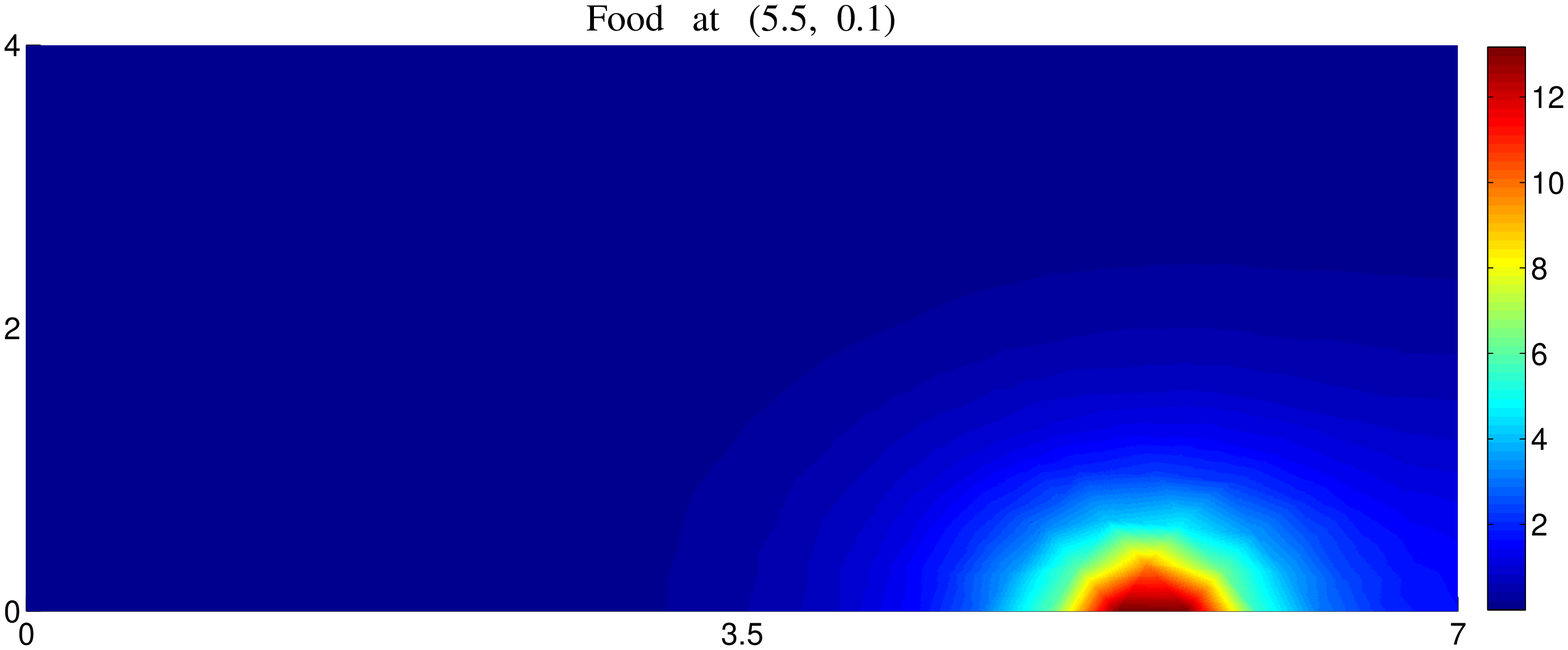} 
 \caption{The color map of potential functions of food scent  $U$ of  \eqref{eq2} in the domain $\Omega$,  where $c=0.1, a=0.2$ and $f$ is defined by \eqref{foodresource}.} 
  \label{ConfigI_food}
 \end{center}
\end{figure}

We put a food resource at a small circle of radius $0.04$ and center being either   $C_1=(1.5;0.1)$ or $C_2=(5.5;0.1)$.  
 More precisely, the function $f$ of food resource in \eqref{eq2} has the form:
\begin{equation}    \label{foodresource}
\begin{aligned}
f(x)=\begin{cases}
50   &  \hspace{1cm} \text{if  }  x\in \mathcal B=\{y \in \mathbb R^2: \|y- C_i\|\leqslant 0.04\}, (i=1 \text{ or } 2),\\
0 &  \hspace{1cm} \text{else}.
\end{cases}
\end{aligned}
\end{equation} 
In addition, set $c=0.1$ and $a=0.2.$
 Then, the elliptic equation  \eqref{eq2} can be numerically solved. Its solution, i.e., the scent function $U,$ is illustrated in Figure \ref{ConfigI_food}. Note that by the Neumann condition,  the scent of food  cannot pass through the boundary $\partial \Omega$. 

We  set initial values and parameters for the system  \eqref{eq4} as follows. All initial positions of $N$-fish are taken randomly in the rectangle domain $[0,2]\times [3.5,4]$, meanwhile all initial velocities are null vectors. Furthermore, $\alpha=\beta=\gamma=1$, $p=P=3$, $q=Q=5$, $r=0.1$, $R=0.2,$ $k=0.5,$ and $\sigma=0.001.$

We introduce a parameter $\|v\|_{\max}$ to restrict speed of fish. If the magnitude of $v_i$ exceeds $\|v\|_{\max}$, our program would reset $v_i$ to a vector whose magnitude is $\|v\|_{\max}$ and direction remains. That is  
\begin{equation*}
v_i(t)=
\begin{cases}
v_i(t)\qquad &\text{\rm if} \quad \|v_i(t)\|\leqslant \|v\|_{\max},\\
\frac{v_i(t)}{\|v_i(t)\|}\|v\|_{\max}\qquad &\text{\rm otherwise}.
\end{cases}
\end{equation*}
This is reasonable because every species owns a maximal speed. 
In our simulations,   $\|v\|_{\max}=0.8$.

We are now ready to perform simulations. We fix all the parameters above except the number $N$ of fish in school, which varies from 2 to 20. (The case $N=1$ is not considered because we want to investigate  behavior of fish school with mutual interaction.) In addition, we perform 100 trials for each $N$.

We say a fish school  succeeds in foraging for the food at time $T=120$ if the distance from the center of the school (i.e. $\bar x=\frac{1}{n} \sum_{i=1}^n x_i$) at that time to the center of food resource $ C_i \,(i=1 \text{ or } 2)$  is less than 1. 

\begin{figure}[H]
\begin{center}
\includegraphics[width=6cm, height=5cm]{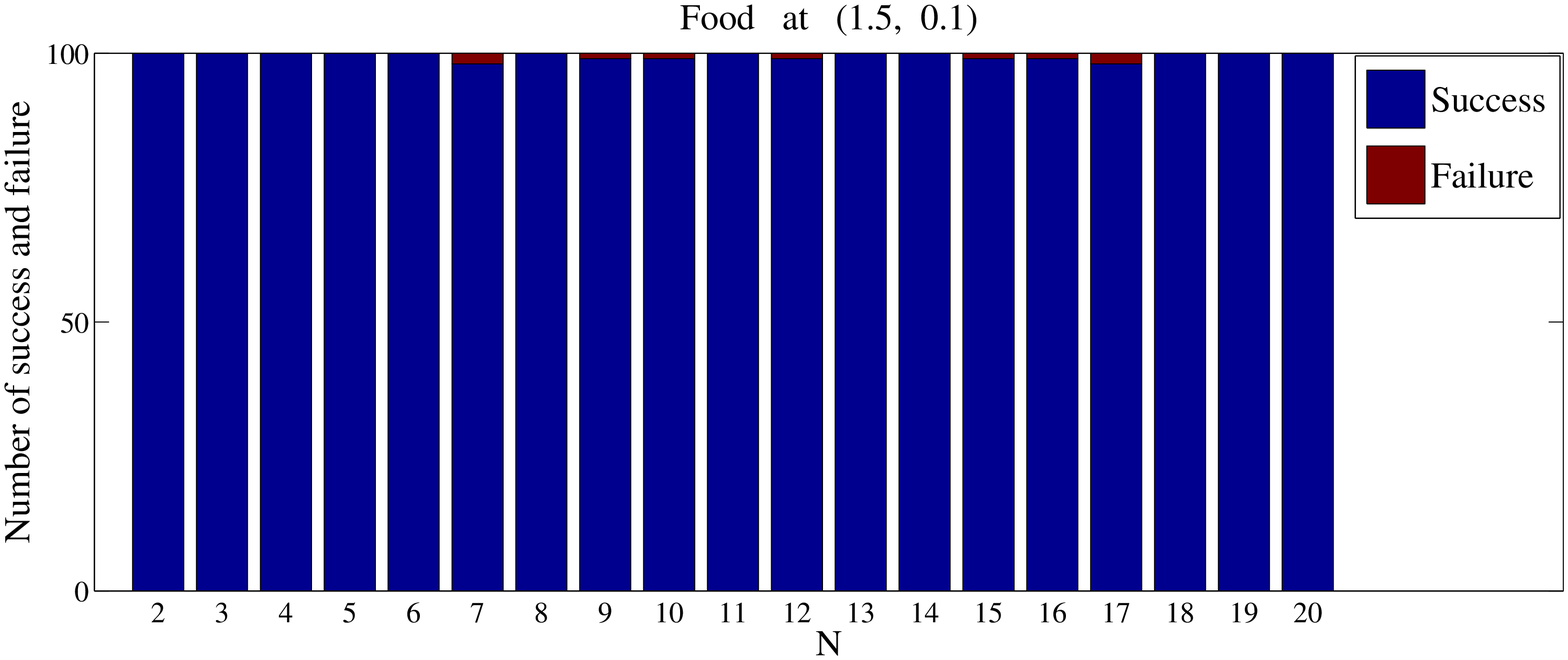} 
\includegraphics[width=6cm, height=5cm]{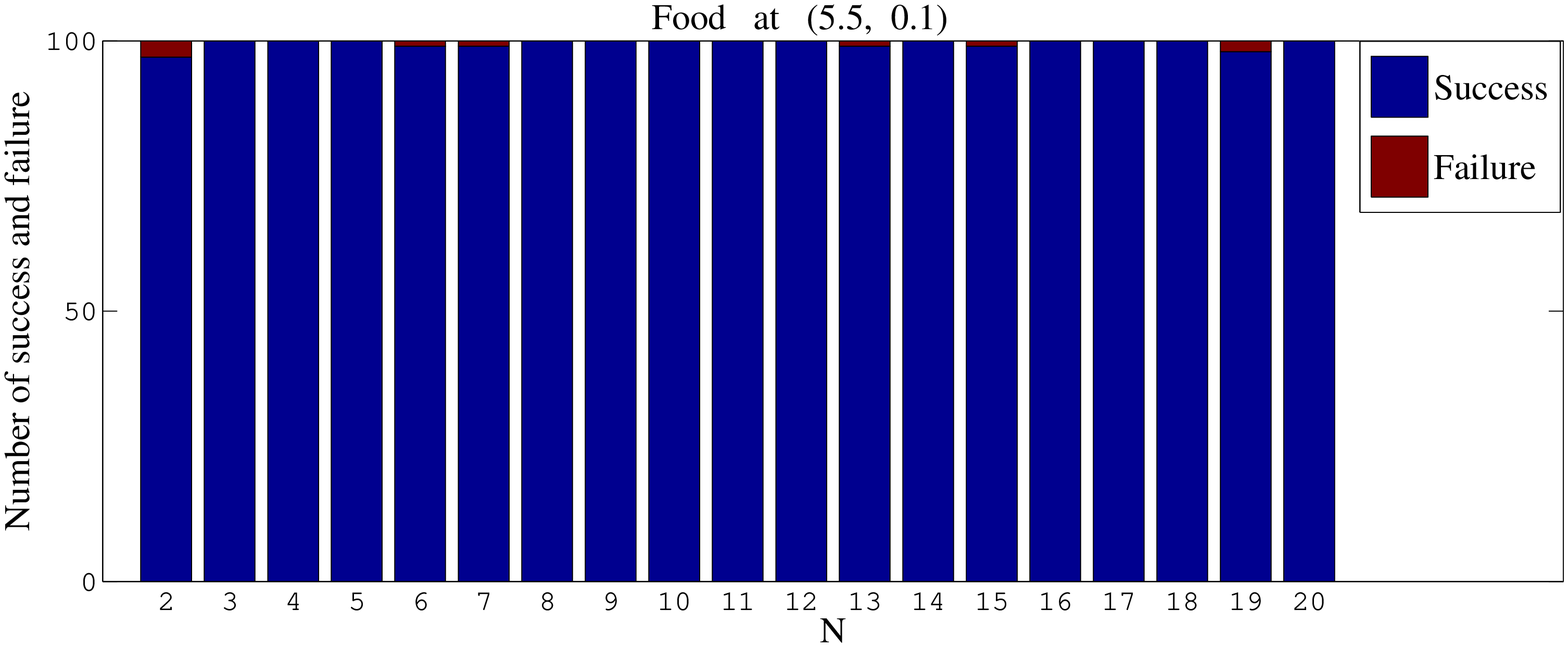} 
 \caption{Numbers of  success and failure in 100 trials. Number of fish in school varies from 2 to 20.} 
%
%
%
  \label{ConfigI_column}
 \end{center}
\end{figure}

\begin{figure}[H]
\begin{center}
\includegraphics[width=6cm, height=5cm]{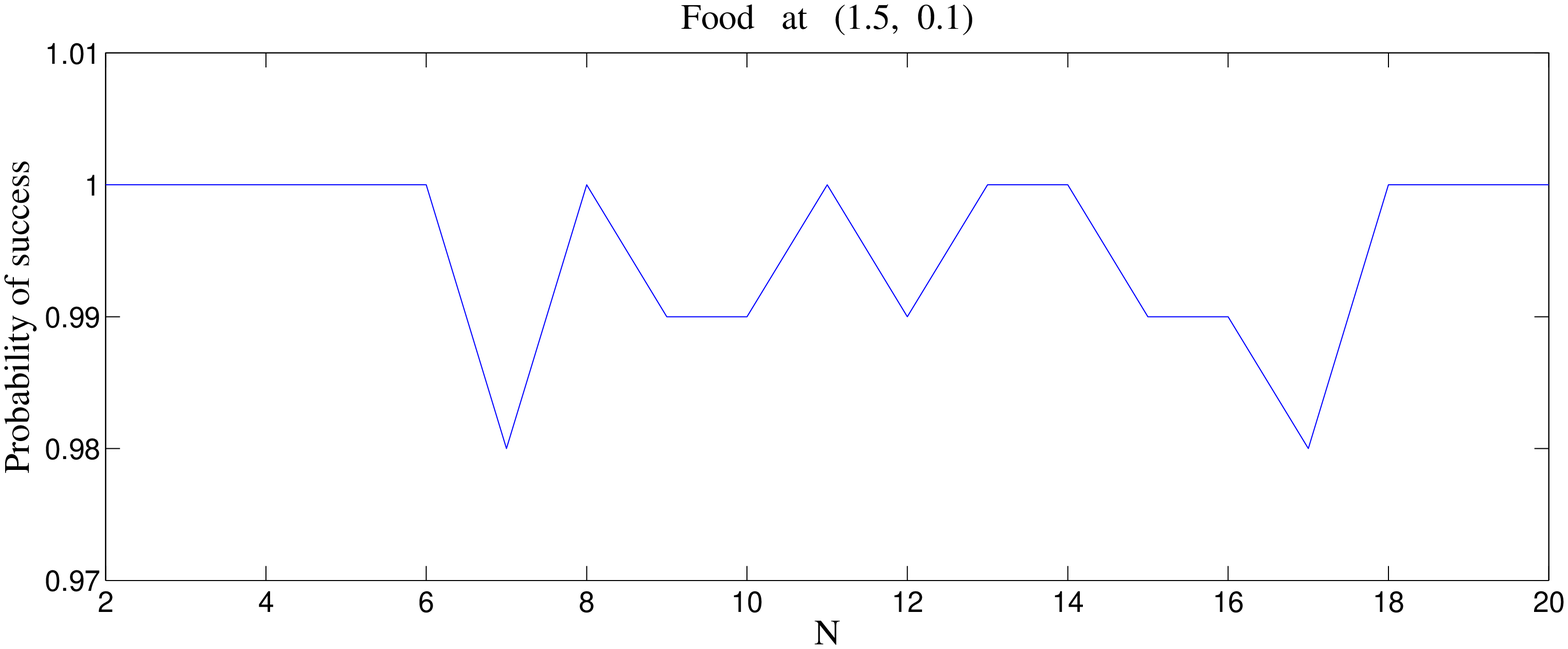} 
\includegraphics[width=6cm, height=5cm]{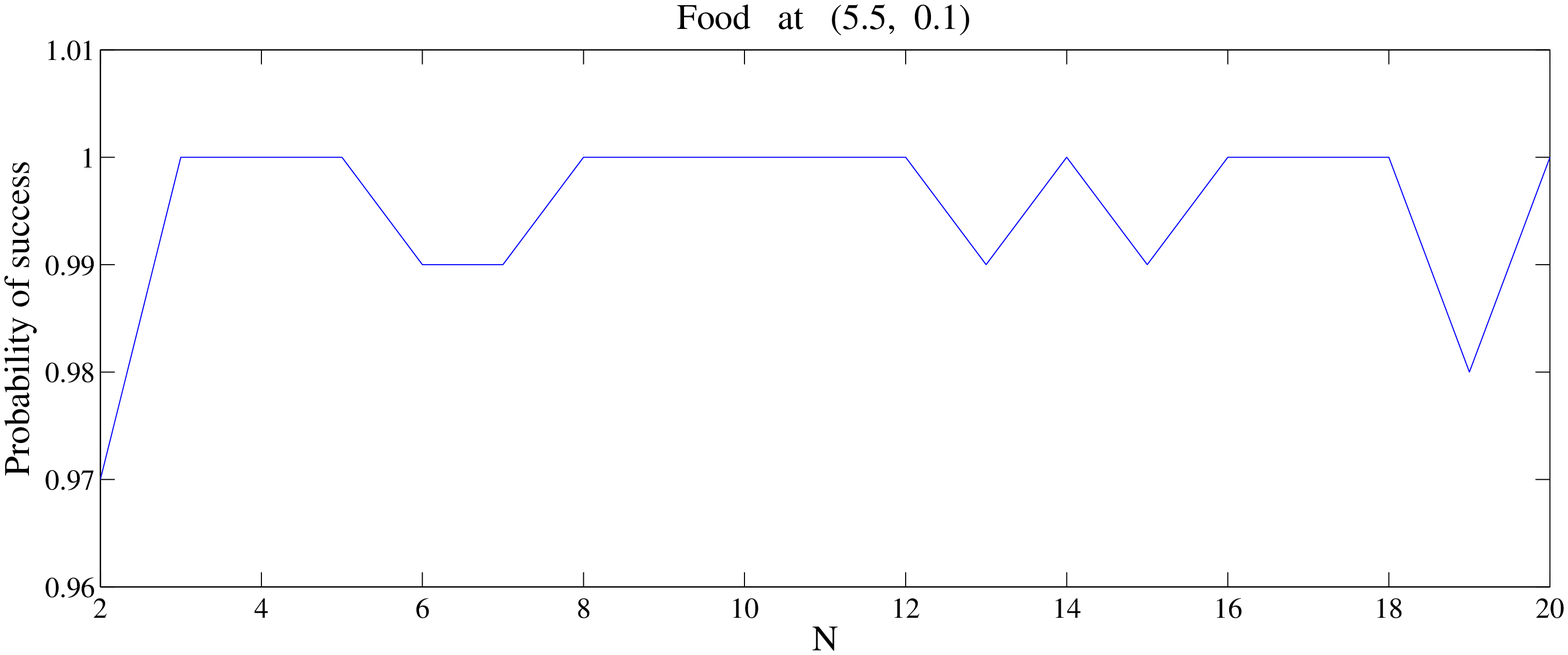} 
 \caption{Probability of success when food is at $(1.5;0.1)$  and $(5.5;0.1).$ Number of fish in school varies from 2 to 20. The probability is calculated based on the data in Figure \ref{ConfigI_column}.}
  \label{ConfigI_success}
 \end{center}
\end{figure}

Figure \ref{ConfigI_column} shows numbers of  success and failure in 100 trials. Graph of the probability function  of foraging success with respect to $N$ is illustrated in Figure \ref{ConfigI_success}. Based on these results, we can say that fish schools almost surely reach the food resource in the domain $\Omega$.

\begin{figure}[H]
\begin{center}
\includegraphics[width=6cm, height=3.5cm]{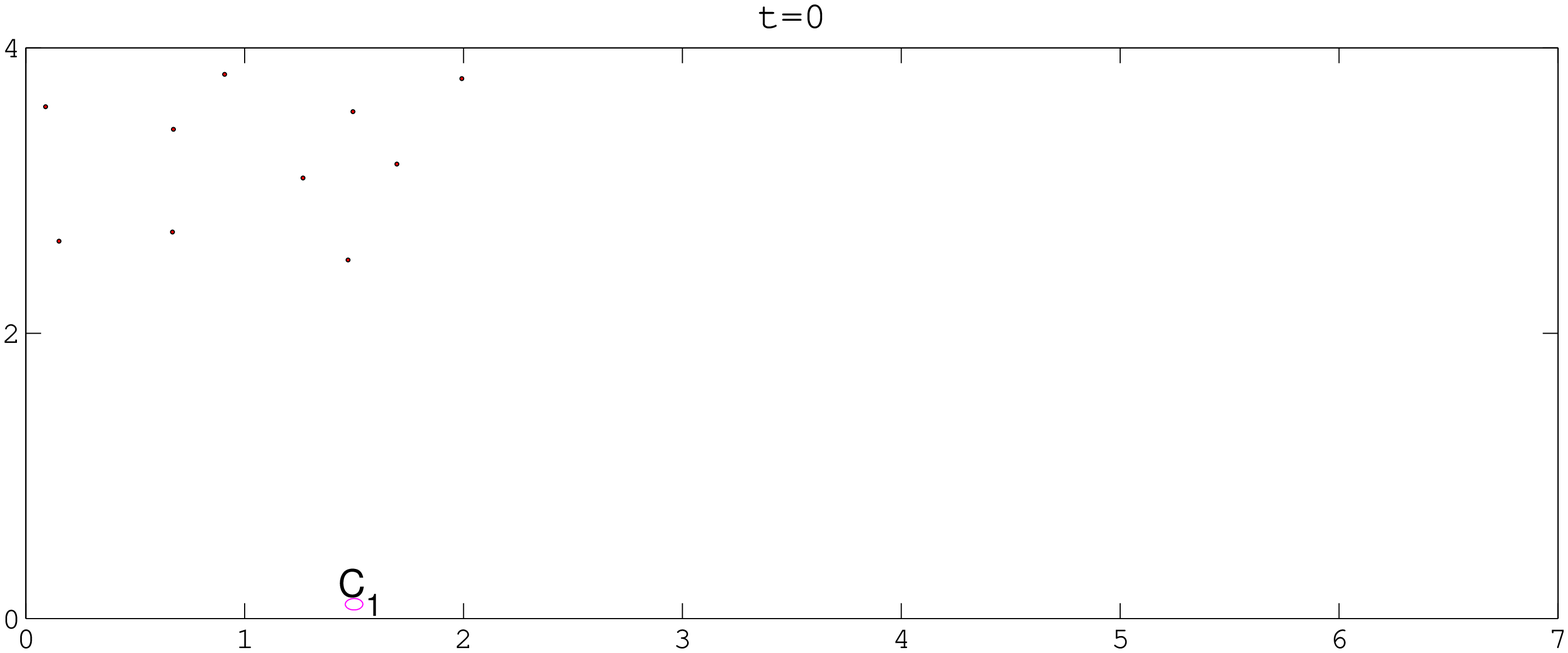} 
\includegraphics[width=6cm, height=3.5cm]{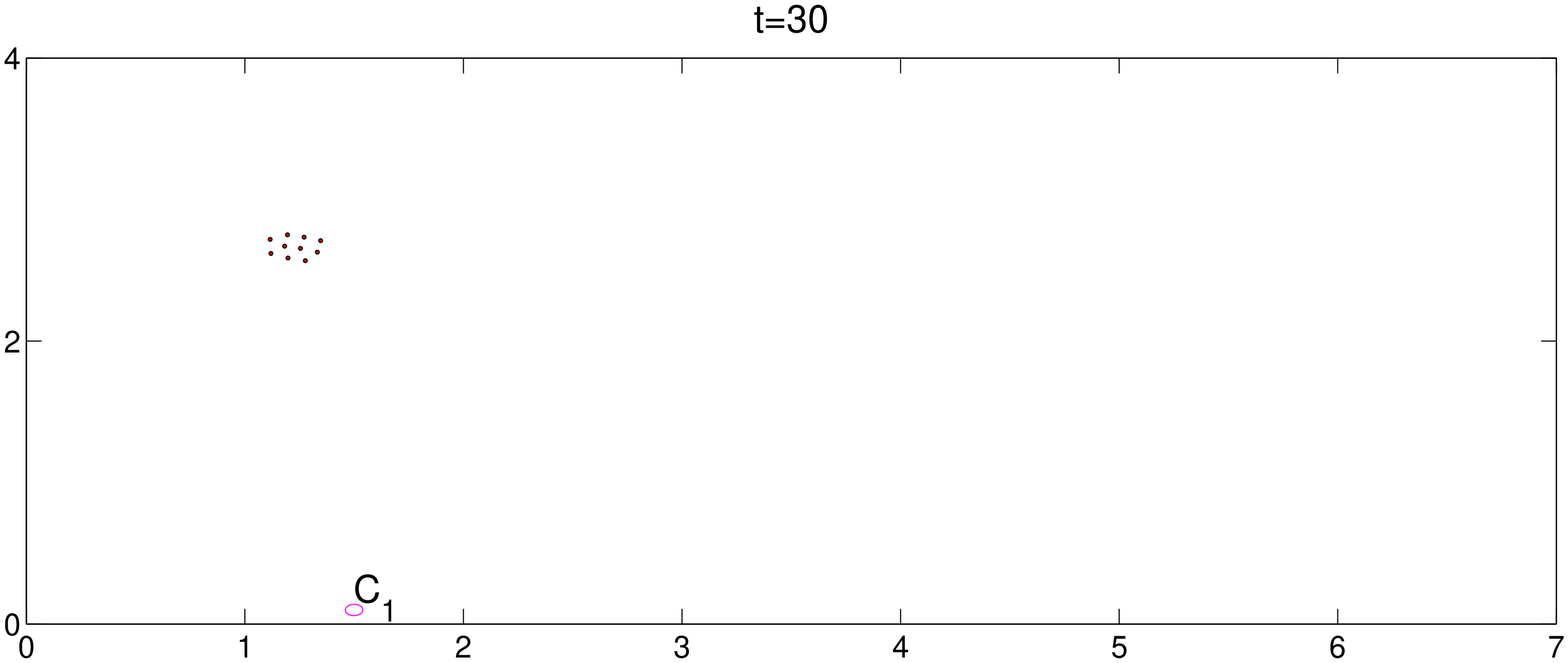} 
\includegraphics[width=6cm, height=3.5cm]{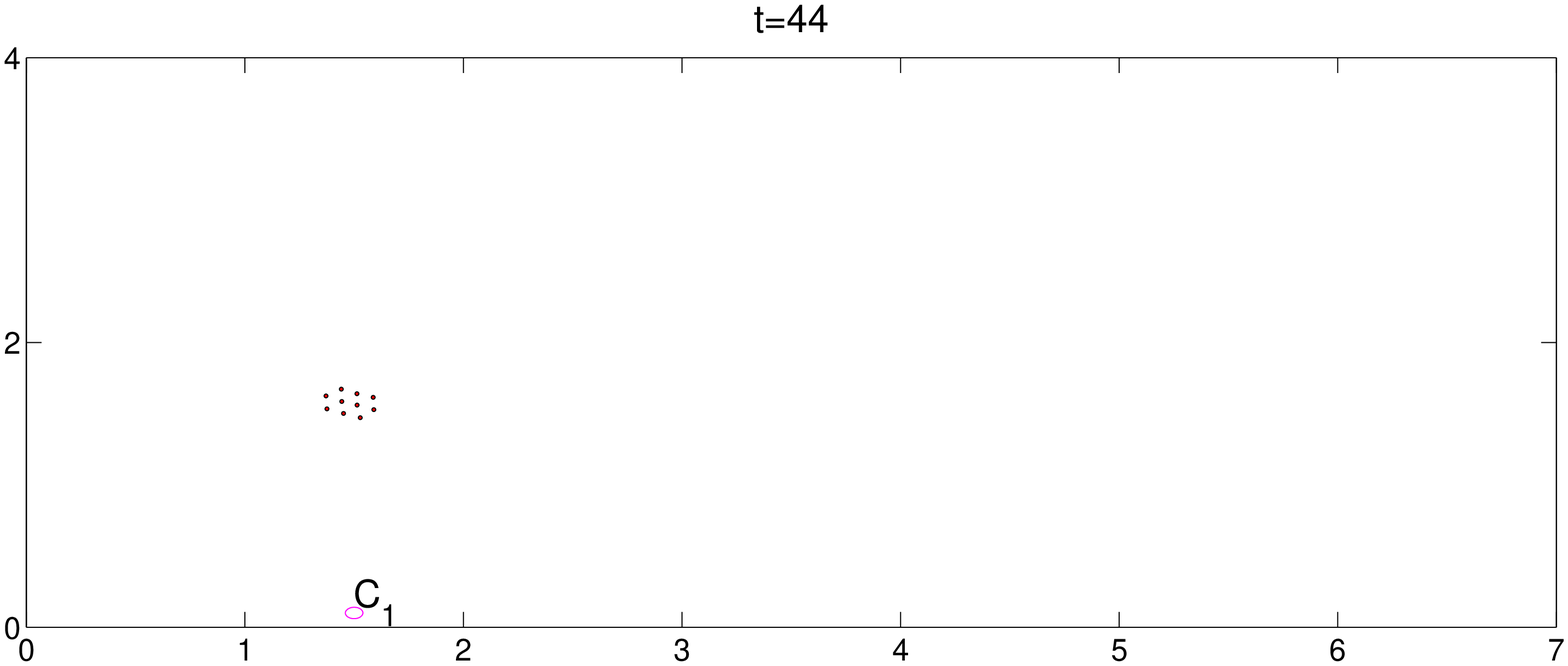} 
\includegraphics[width=6cm, height=3.5cm]{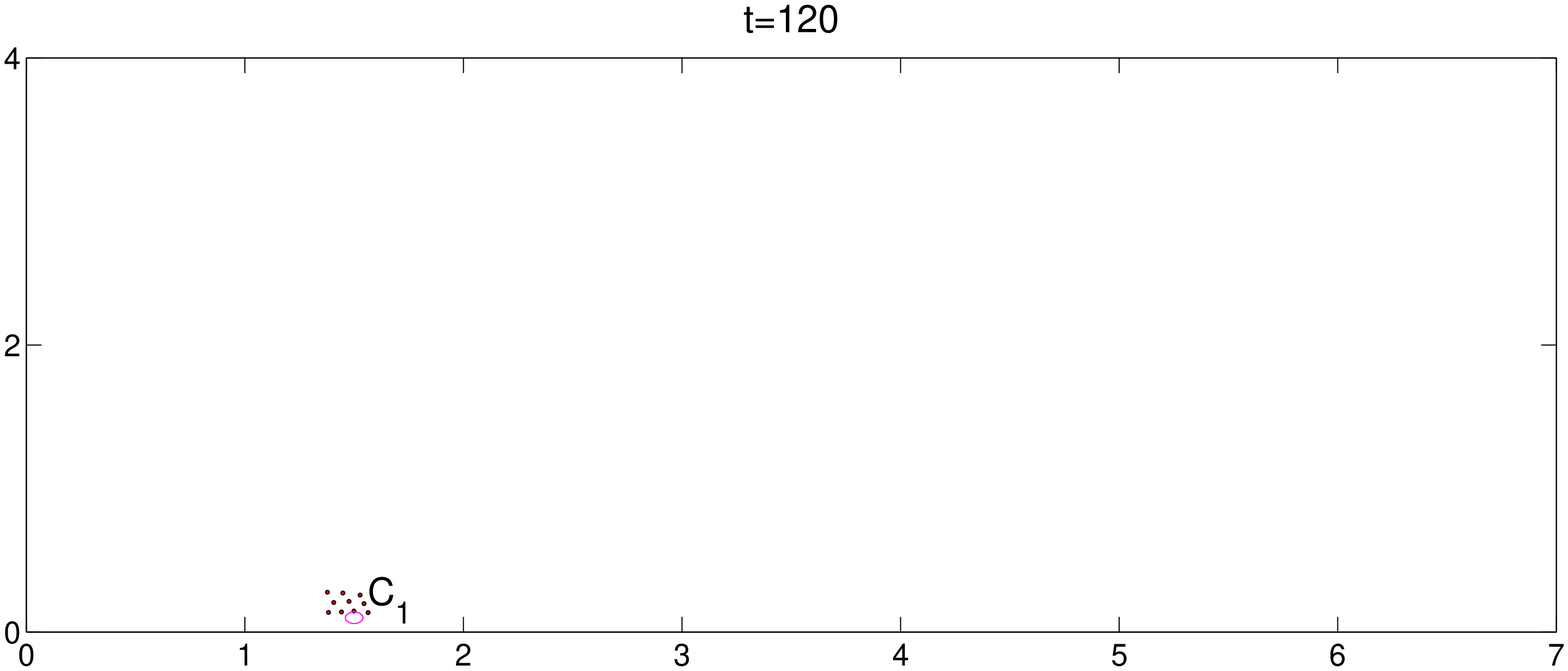} 
 \caption{A pattern of collective foraging in free but limited space. Positions of all ten fish in school at four instants are plotted. The food resource is at $C_1=(1.5;0.1).$ The school reaches to the food resource while keeping its school structure during the process of foraging.} 
  \label{ConfigI_left_position}
\vspace{0.4cm}
\includegraphics[width=6cm, height=3.5cm]{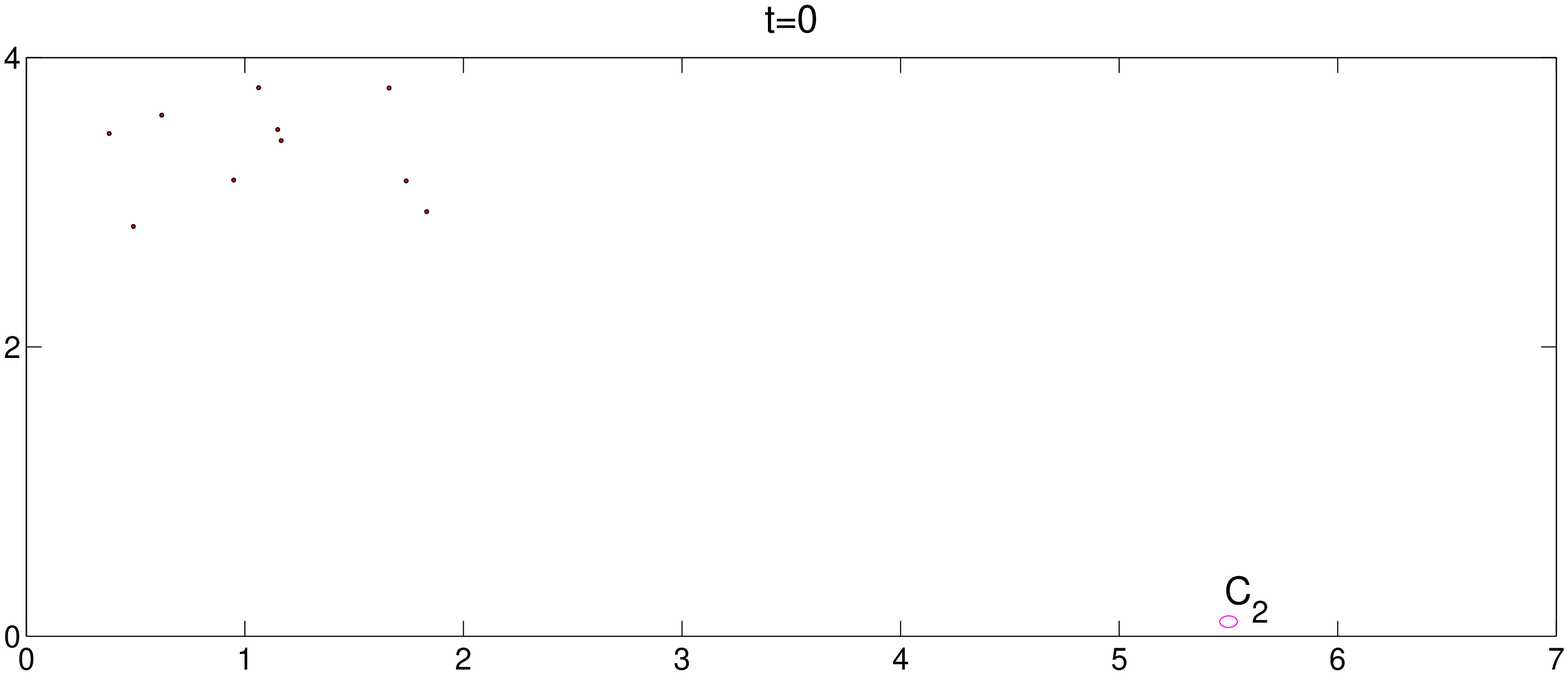} 
\includegraphics[width=6cm, height=3.5cm]{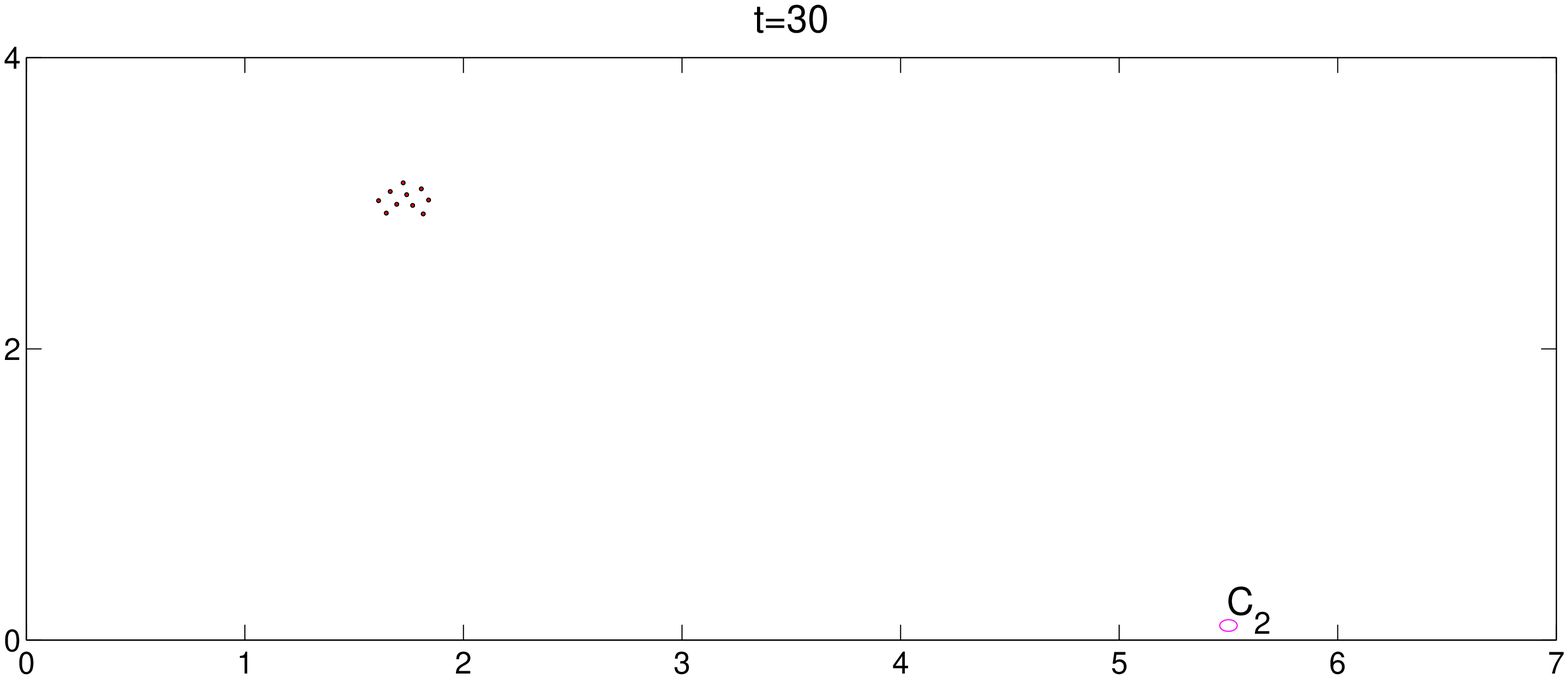} 
\includegraphics[width=6cm, height=3.5cm]{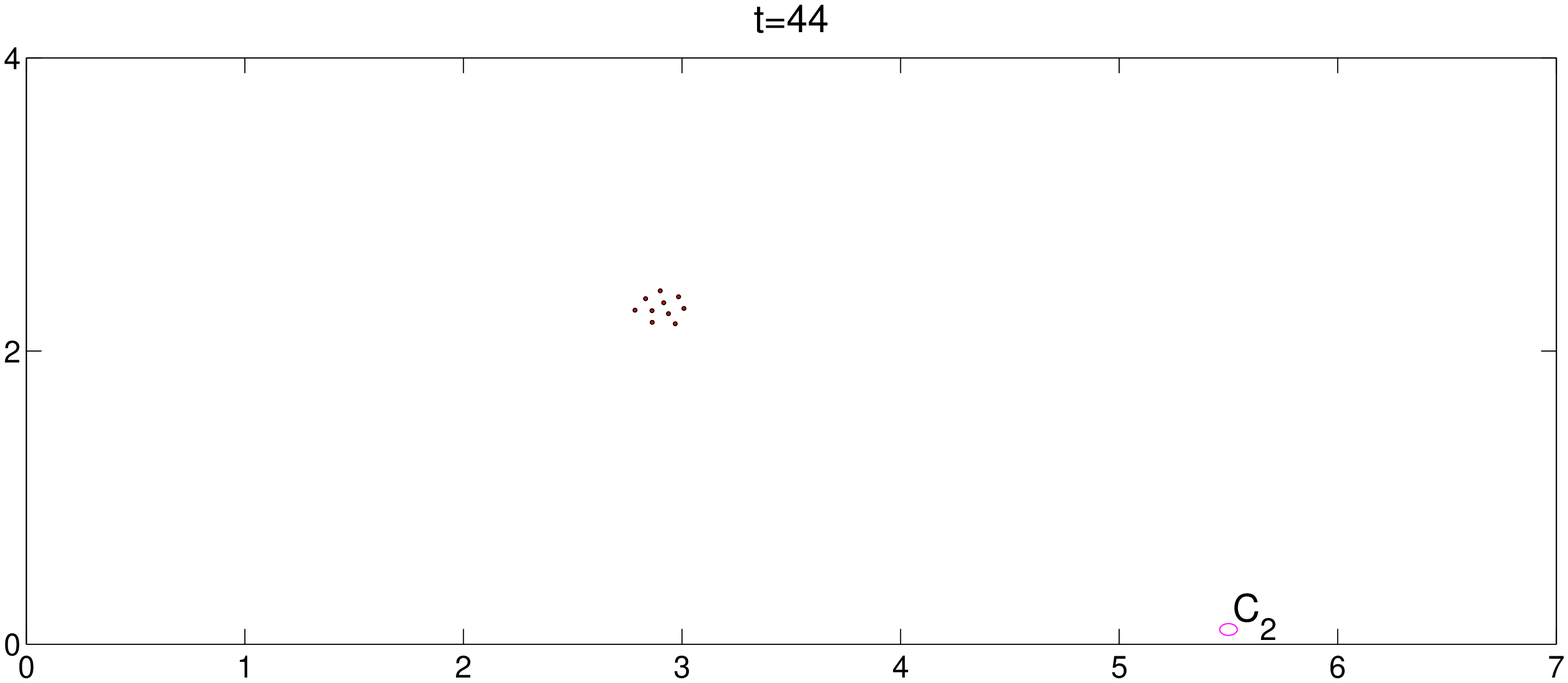} 
\includegraphics[width=6cm, height=3.5cm]{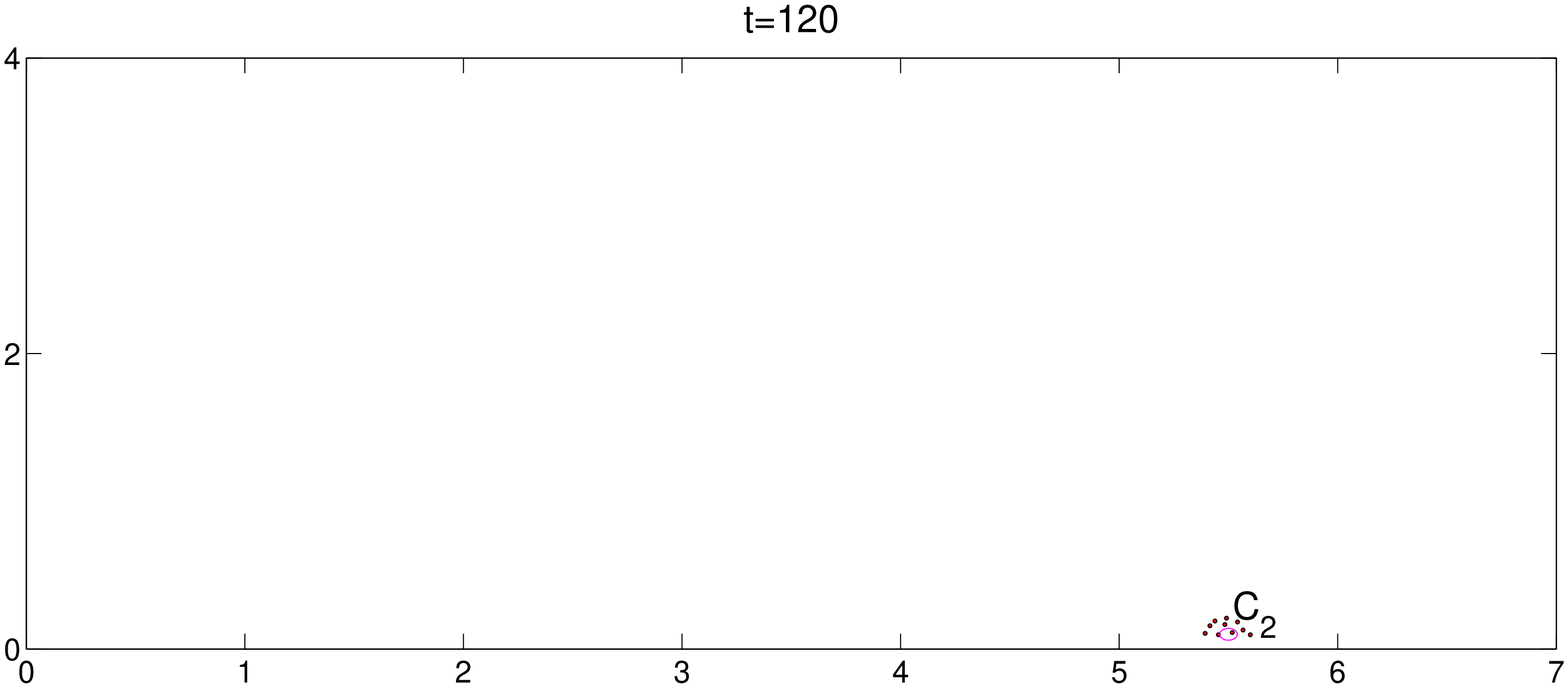} 
 \caption{A  pattern of collective foraging in free but limited space. It is similar to the pattern in Figure \ref{ConfigI_left_position} but the food resource here is at $C_2=(5.5;0.1).$} 
  \label{ConfigI_right_position}
 \end{center}
\end{figure}

Figures \ref{ConfigI_left_position} and \ref{ConfigI_right_position} show some patterns of collective foraging at time $t=0, 30, 44,$ and $120$ by using one of the above 100 trials.  It is seen that fish always keep school structure during the process of foraging.

\subsection{Configuration II}

Let us consider the system  \eqref{eq4} in a domain $\Omega=S\setminus Ob_1$, where
$$S=[0, 4] \times [0, 4],$$
and $Ob_1$ is an obstacle, say a long thin rectangle put inside $S$:   
$$Ob_1=[2, 2.5]\times[2.5,4].$$

\begin{figure}[H]
\begin{center}
\includegraphics[width=10cm, height=6cm]{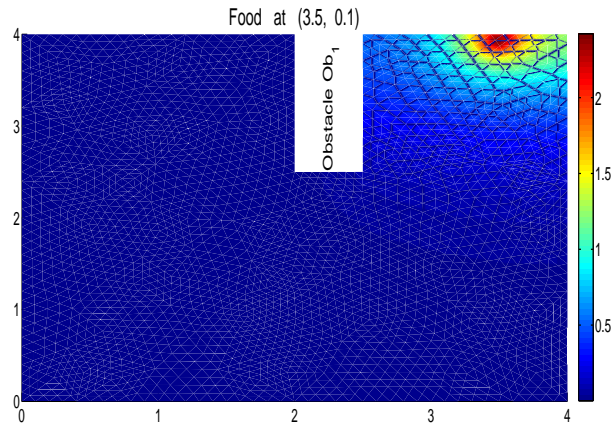} 
 \caption{The color map of potential functions of food scent  $U$ of  \eqref{eq2} in the domain $\Omega=S\setminus Ob_1$,  where $c=0.1, a=0.2$ and $f$ is defined by \eqref{foodresource} (with $i=3$).}
  \label{up_odour}
 \end{center}
\end{figure}

The food resource is put at a small circle of radius $0.04$ and center  $C_3=(3.5;0.1)$. The function $f$ of food resource (with $i=3$ in \eqref{foodresource}) and the parameters $a$ and $c$ in \eqref{eq2} are the same  as in the configuration I. The scent function $U$ in  \eqref{eq2} can be then numerically solved in $\Omega$. The color map of  $U$ is illustrated in Figure \ref{up_odour}. It is similar to  the configuration I that the scent of food  cannot pass through the boundary $\partial \Omega$.

Other parameters in \eqref{eq4} as well as the parameter $\|v\|_{\max}$ are the same as in the configuration I except for the followings. Initial positions of $N$-fish are taken randomly in the rectangle domain $[1,2]\times [3.5,4]$. Allotted time is $T=60$, and the sensitivity constant is $k=2.$  

In our simulations, we do 200 trials to the system \eqref{eq4} for each $N$ (from 2 to 20).  In view of the configuration I,   schools of fish reach food resource almost surely in free-obstacle domain with $k=0.5$. Therefore, once a fish school in this configuration ($k=2>0.5$) has moved to the right-hand side of  $Ob_1,$ the school would certainly succeeds in foraging for food. Thus, we 
 classify the state of a fish school at the allotted time $T=60$ into 2 states:

{\bf State I (Failure)}:  Some fish in the school are on the left-hand side of the right wall of $Ob_1.$ A mathematical expression for this state is that  
$$\min_{i=1,2\dots N} x_i^1(T) \leq 2.5,$$ 
where $x_i^1$ is the first component of vector $x_i$.

{\bf State II (Success)}:    The whole school are on  the right-hand side of  $Ob_1$, i.e., 
$$\min_{i=1,2\dots N} x_i^1(T) >2.5.$$

\begin{figure}[H]
\begin{center}
\includegraphics[width=12cm, height=6cm]{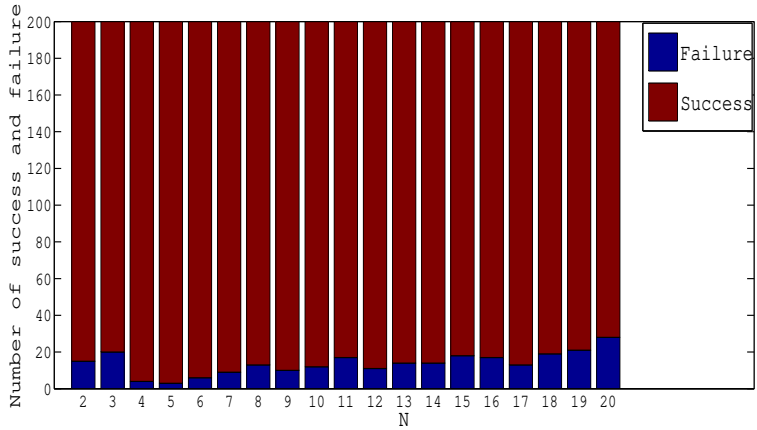} 
 \caption{Numbers of  success and failure in 200 trials at the allotted time $T=60$. Number of fish in school varies from 2 to 20.} 
  \label{k2_column_Config3}
 \end{center}
\end{figure}

Numbers of success and failure in 200 trials are shown in Figure \ref{k2_column_Config3}.

\begin{figure}[H]
\begin{center}
\includegraphics[width=12cm, height=6cm]{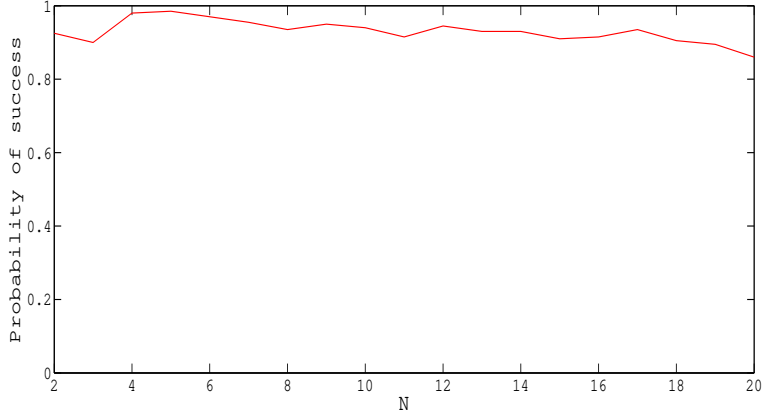} 
 \caption{Probability function of foraging success with respect to $N$.  The probabilities are calculated based on the data in Figure \ref{k2_column_Config3}, where the food is at $C_3=(3.5, 4)$.} 
  \label{Config3_success}
 \end{center}
\end{figure}

Based on the data in Figure \ref{k2_column_Config3}, graph of the probability function  of  success with respect to $N$ is illustrated in Figure \ref{Config3_success}.  It is then seen that on one hand, when size of fish school increases up to some optimal value, so does the probability of foraging success. 
In our simulations, this optimal value is $N=5.$ 

The fact may be explained as follows. Since  fish take initial positions in a small rectangle  in the  left-hand side of the obstacle, where the scent of food is weak,  one fish itself senses the scent weakly. If there are more fish in school, some  can be at good positions to smell food better. These fish have tendency to swim towards the food resource along the gradient of  scent. Their neighbor fishes would then follow them due to social interactions in the model (the attraction term). This movement results in the success of the whole school in foraging for food. Another possible reason can be that the total effect of food scent on a big school of fish is stronger than that on a smaller one. 

On the other hand, Figures \ref{k2_column_Config3}-\ref{Config3_success} reveal that  the probability of success  trends to decrease as $N$ exceeds the optimal value. We may explain this fact as follows. 
In \cite{Nguyen2014}, we showed that number of {\it connected components} of fish school decreases as school size increases due to the model \eqref{eq0}. (A connected component is a set of particles in which for any particle there exists a nearby particle such that the distance between the two particles is less than a given small constant.) Furthermore, it is shown that all fish in school are connected if school size is greater than some value.  
(On the other words, we can say {\it school cohesiveness},  that is defined in \cite{LinhTonYagi} as the ability of group of fish to form and maintain ``connections between fish"  against noise, increase as school size increases.) 
 Therefore, if only a few fish pass the obstacle, by strong social connections (interactions) between fish due to the attraction force in the model, the crowd of other fish would probably pull them to ``previous positions", i.e., the left-hand side of the obstacle, or at least slows down movement of these fish to the food resource. In other words, the larger the fish school the slower the movement to the food source. Consequently,  the probability of success at the allotted time decreases.

\begin{figure}[H]
\begin{center}
\includegraphics[width=6cm, height=4cm]{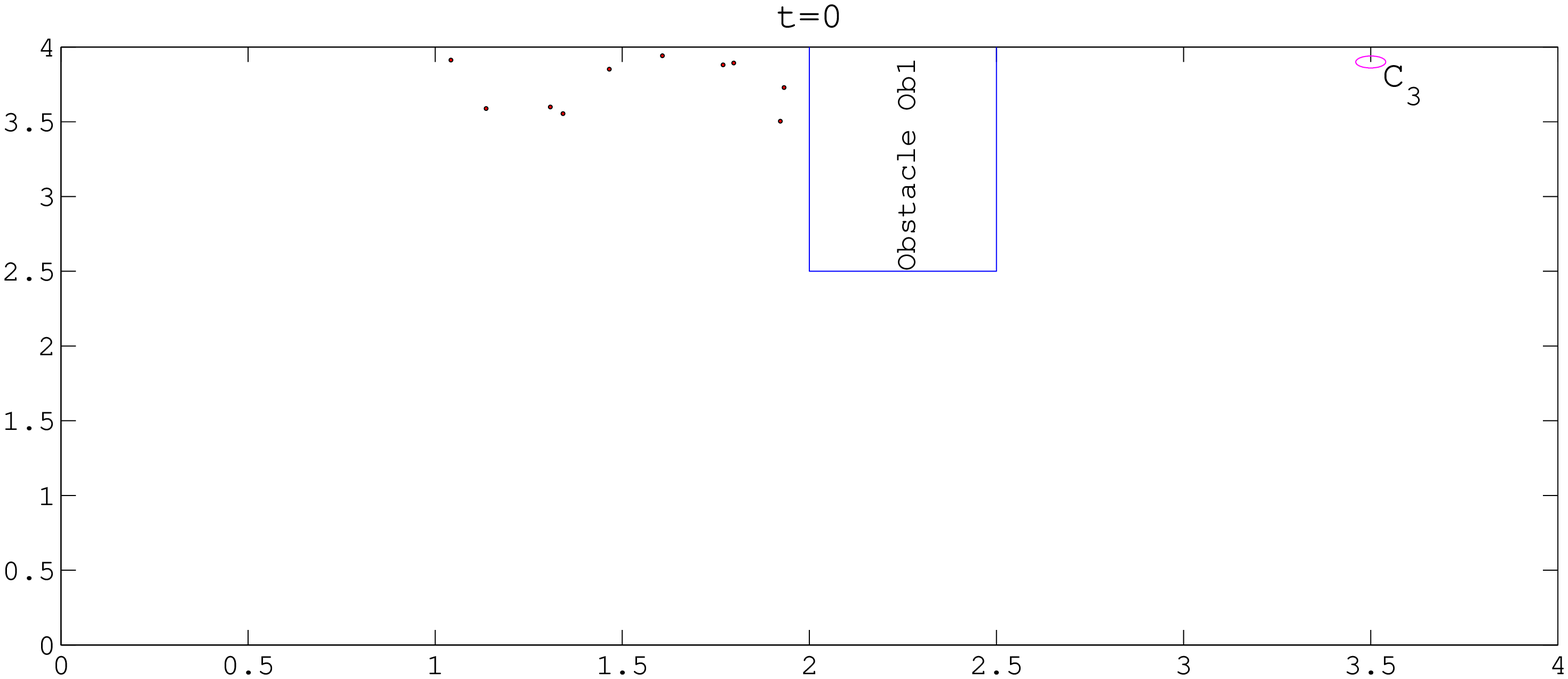} 
\includegraphics[width=6cm, height=4cm]{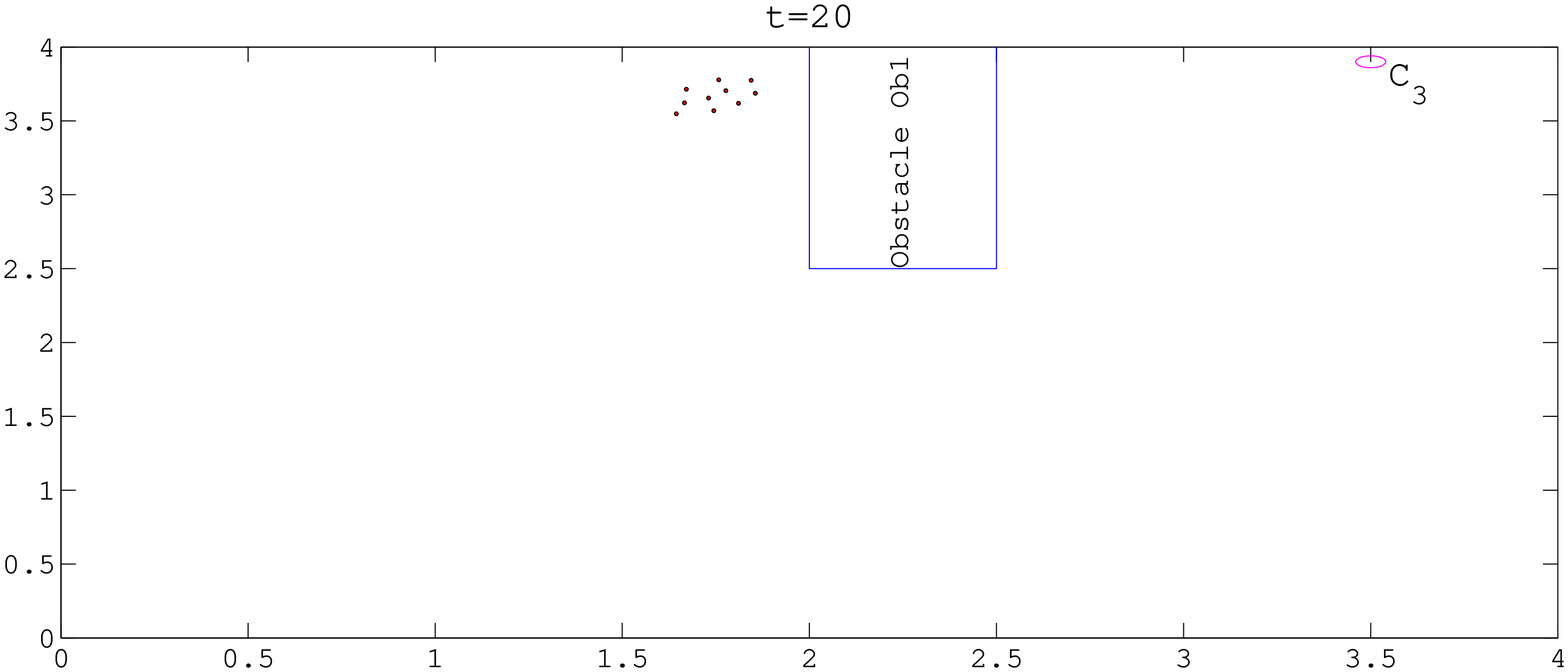} 
\includegraphics[width=6cm, height=4cm]{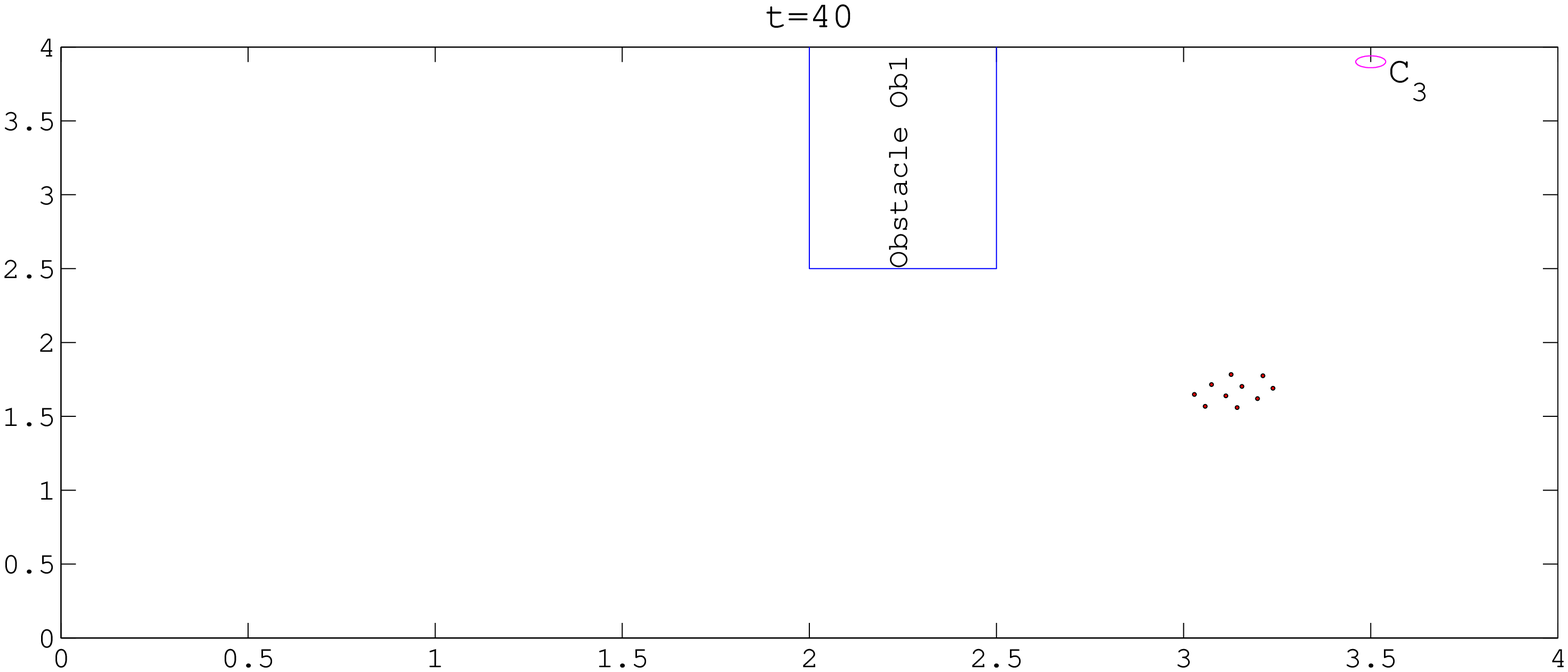} 
\includegraphics[width=6cm, height=4cm]{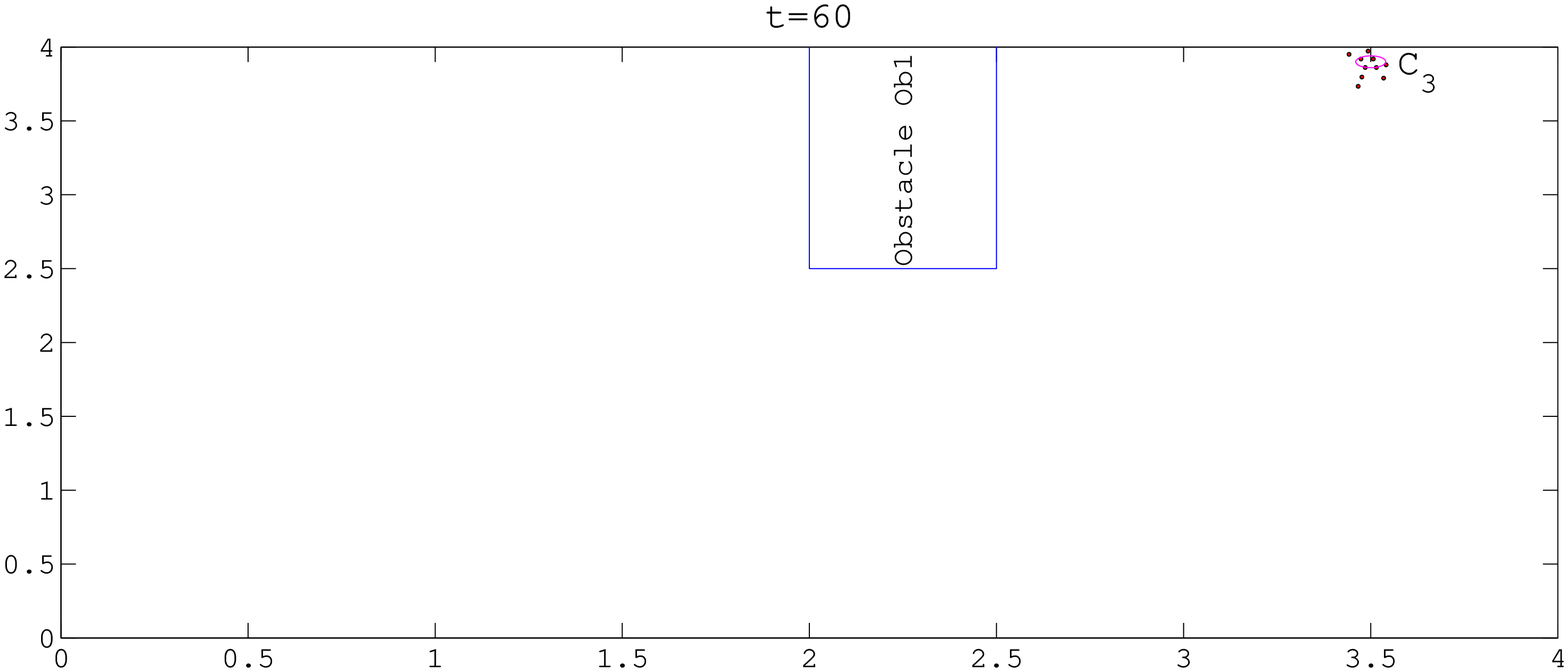} 
 \caption{A pattern of collective foraging. Positions of all ten fish in school are demonstrated at four instants. The school succeeds in foraging for food putting at $C_3=(3.5;0.1)$, and maintains its school structure during the process of foraging.} 
  \label{Config3_position}
 \end{center}
\end{figure}

In Figure \ref{Config3_position}, we give a pattern of collective foraging. Positions of all fish in school at time $t=0, 20, 40,$ and $60$ are plotted by using one of the above 200 trials.

\subsection{Configuration III}

In this configuration, we consider the system  \eqref{eq4} in the  domain $\Omega$ defined by
$$\Omega=[0, 7]\times [0, 4]\setminus (Ob_1 \cup Ob_2),$$
 where $Ob_1$ is the same obstacle as in the configuration II, and   $Ob_2$ is another obstacle, say a long thin rectangle  defined by
$$Ob_2=[4.5, 5]\times [0,1.5].$$

\begin{figure}[H]
\begin{center}
\includegraphics[width=10cm, height=6cm]{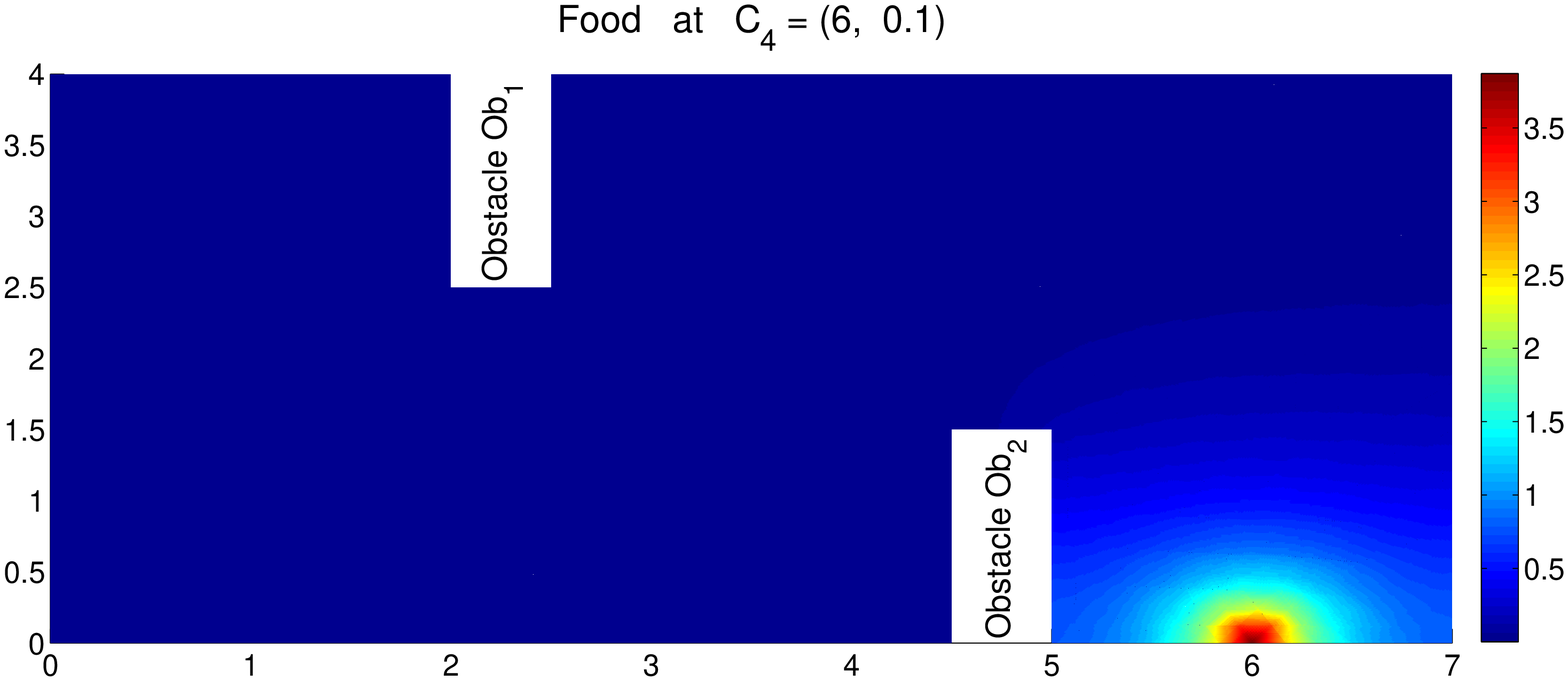} 
 \caption{The color map of potential function of food scent  $U$ of  \eqref{eq2} in the domain $\Omega=[0, 7]\times [0, 4]\setminus (Ob_1 \cup Ob_2)$,  where $c=0.1, a=0.2$ and $f$ is defined by \eqref{foodresource} (with $i=4$).}
  \label{bb2Fig4}
 \end{center}
\end{figure}

Parameters are the same as in the configuration II except position of food and the allotted time.  The food is now at   $C_4=(6, 0.1)$, and the allotted time is $T=200$.  
The numerical solution of \eqref{eq2}   in $\Omega$, i.e., the scent function $U$ is illustrated by its  color map  in Figure \ref{bb2Fig4}. 

In our simulations, we perform 200 trials to the system \eqref{eq4} for each value of $N$ from 2 to 20. As in the configuration II,  we 
 classify the state of a fish school  at the allotted time $T=200$ into 3 states:

{\bf State I (Failure)}:  The whole school is on the left-hand side of the obstacle $Ob_1$. Mathematically, this means that  
$$\max_{i=1,2\dots N} x_i^1(T) <2.$$ 

{\bf State II (Pre-success)}:  Some (or all) fish in the school are in the domain limited by the  left wall of $Ob_1$ and  the right wall of $Ob_2$. Precisely, 
$$2\leq \max_{i=1,2\dots N} x_i^1(T) \quad \text { and }  \quad  \min_{i=1,2\dots N} x_i^1(T) \leq 5.$$ 

{\bf State III (Success)}:    The whole school has moved to  the right-hand side of  $Ob_2$, i.e.
$$\min_{i=1,2\dots N} x_i^1(T) >5.$$

\begin{figure}[H]
\begin{center}
\includegraphics[width=12cm, height=6cm]{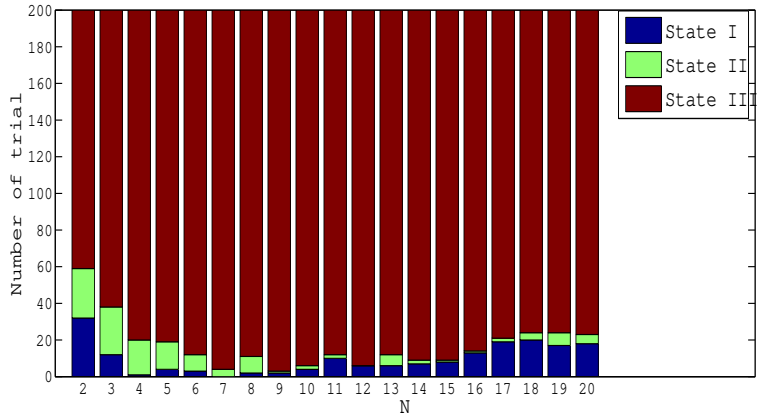} 
 \caption{200 trials to the system \eqref{eq4}  are clarified into the three states I (failure), II (pre-success), and III (success) at the allotted time $T=200.$   Number of fish in school varies from 2 to 20.} 
  \label{k2_column_Config2}
 \end{center}
\end{figure}

Numbers of trial fell into these states are shown in Figure \ref{k2_column_Config2}.

\begin{figure}[H]
\begin{center}
\includegraphics[width=12cm, height=5cm]{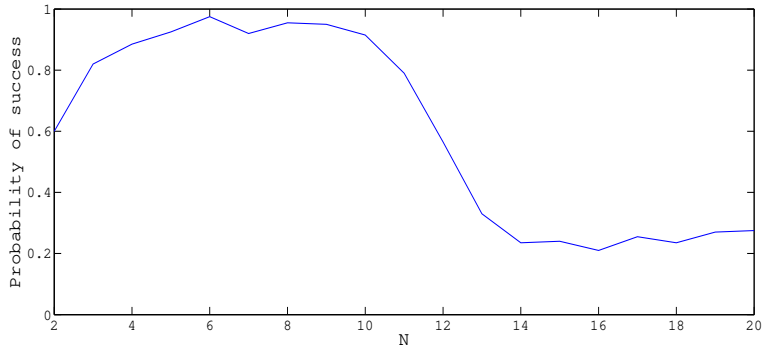} 
 \caption{Probability function of foraging success with respect to $N$.  The probabilities are calculated based on the data in Figure \ref{k2_column_Config2}.} 
  \label{Config2_success}
 \end{center}
\end{figure}

The probability function  of  success with respect to $N$ is then illustrated in Figure \ref{Config2_success}.

From Figures \ref{k2_column_Config2}--\ref{Config2_success}, we obtain the same conclusion as in the configuration II. That is, on one hand, as the size of fish school increases up to an optimal value ($N=9$), so does the probability of foraging success. 
On the other hand,  this probability decreases from this optimal value.

\begin{figure}[H]
\begin{center}
\includegraphics[width=6cm, height=5cm]{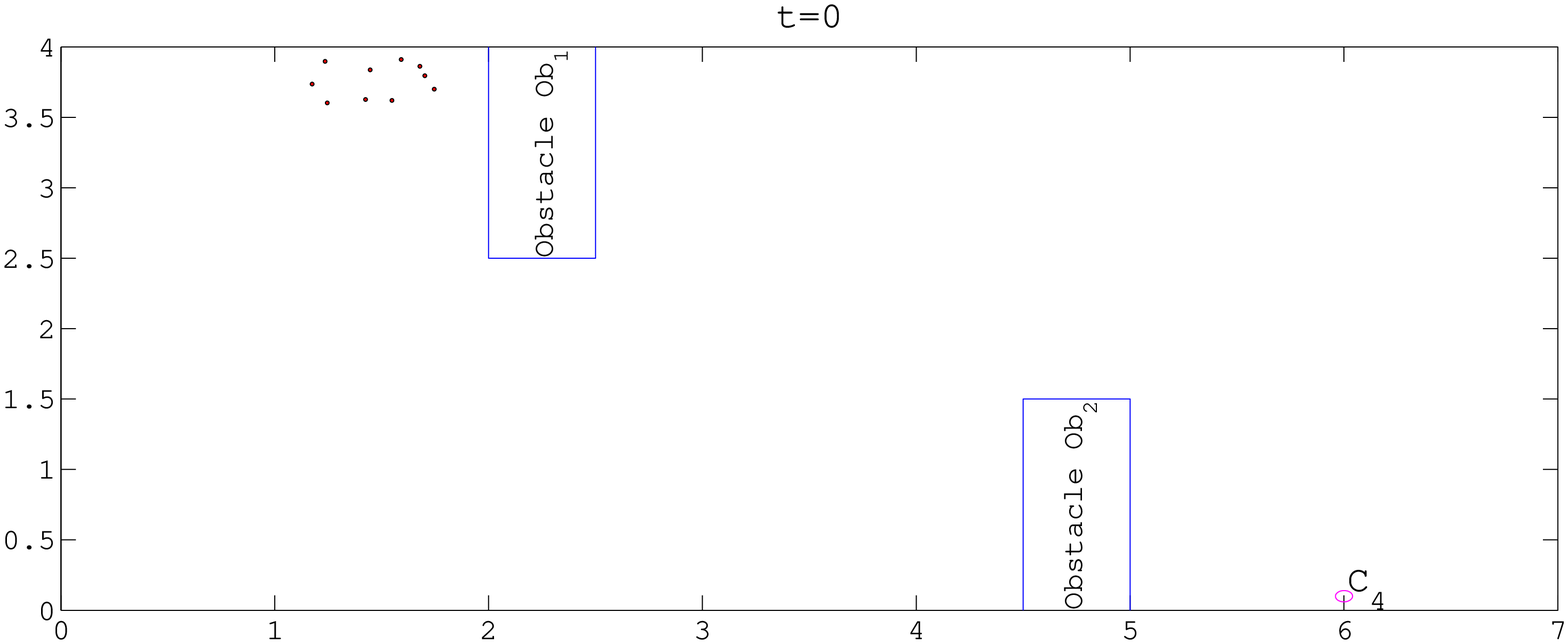} 
\includegraphics[width=6cm, height=5cm]{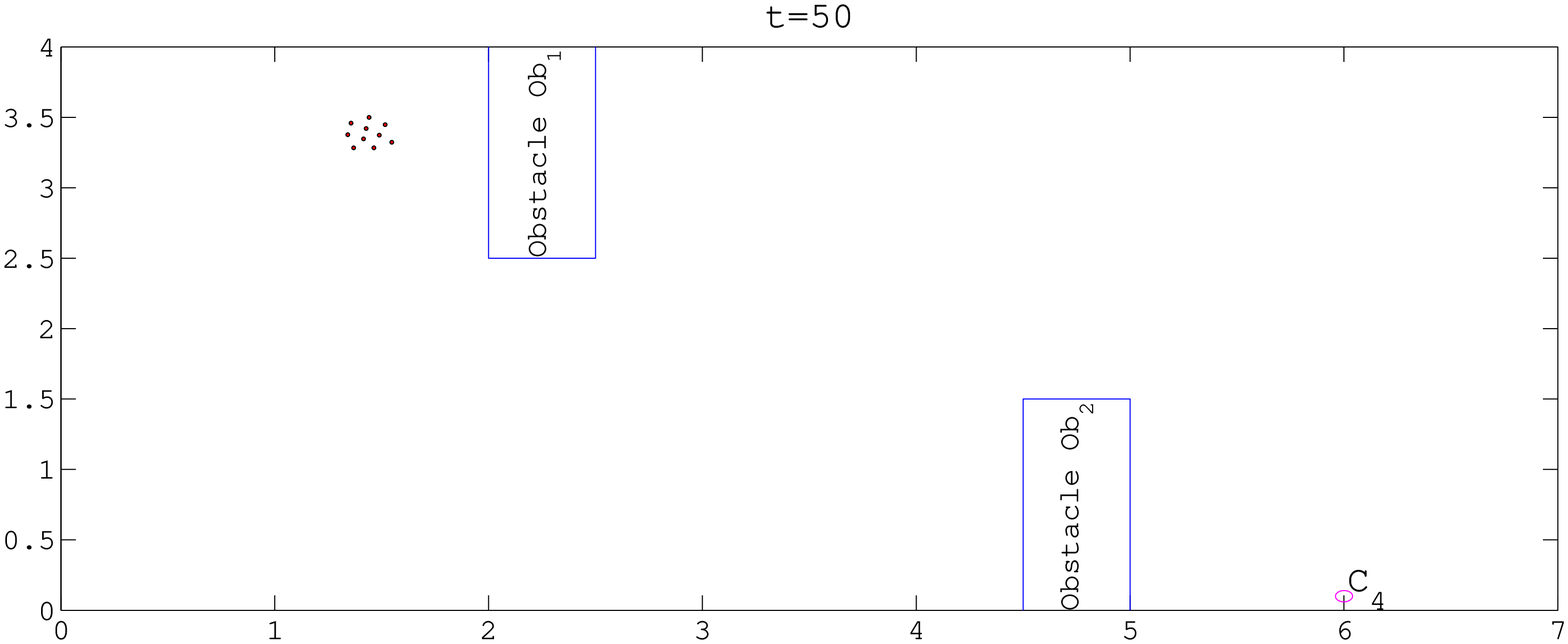} 
\includegraphics[width=6cm, height=5cm]{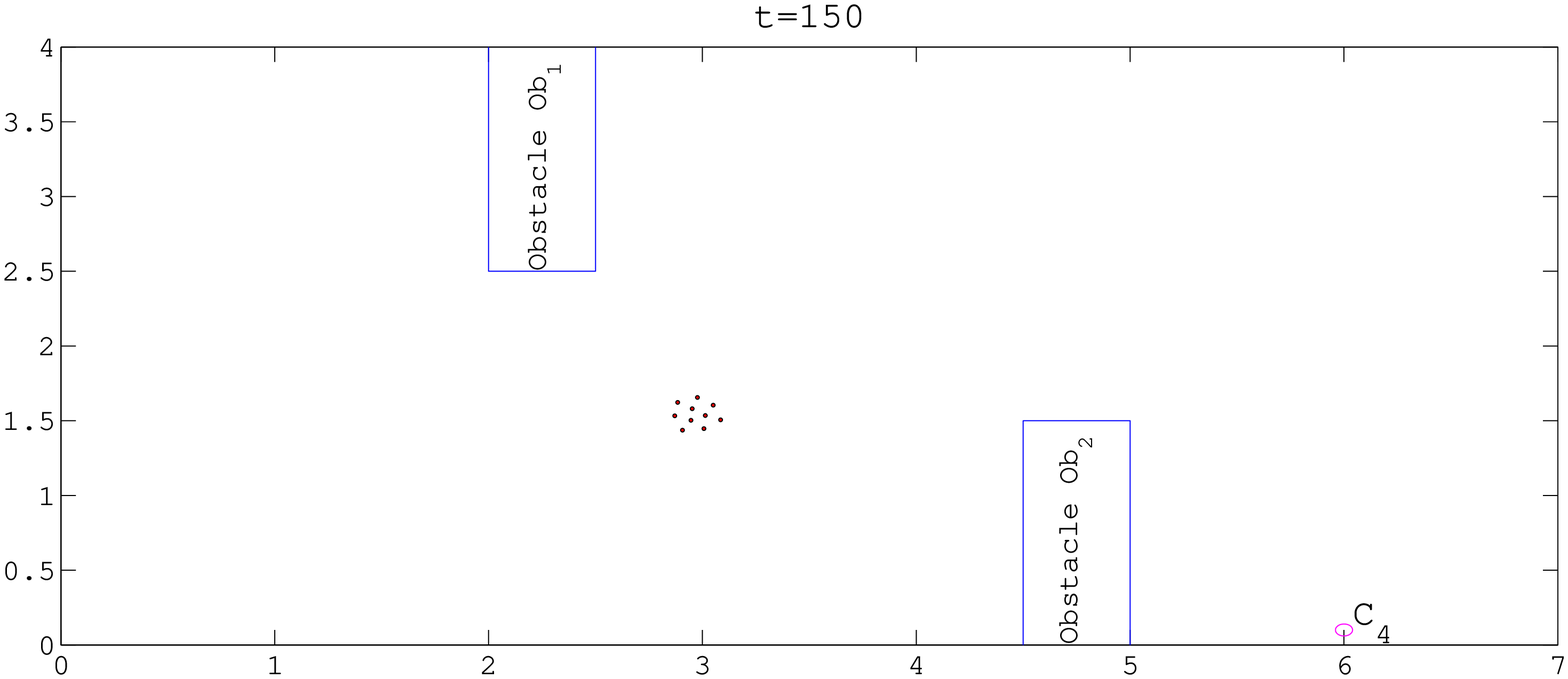} 
\includegraphics[width=6cm, height=5cm]{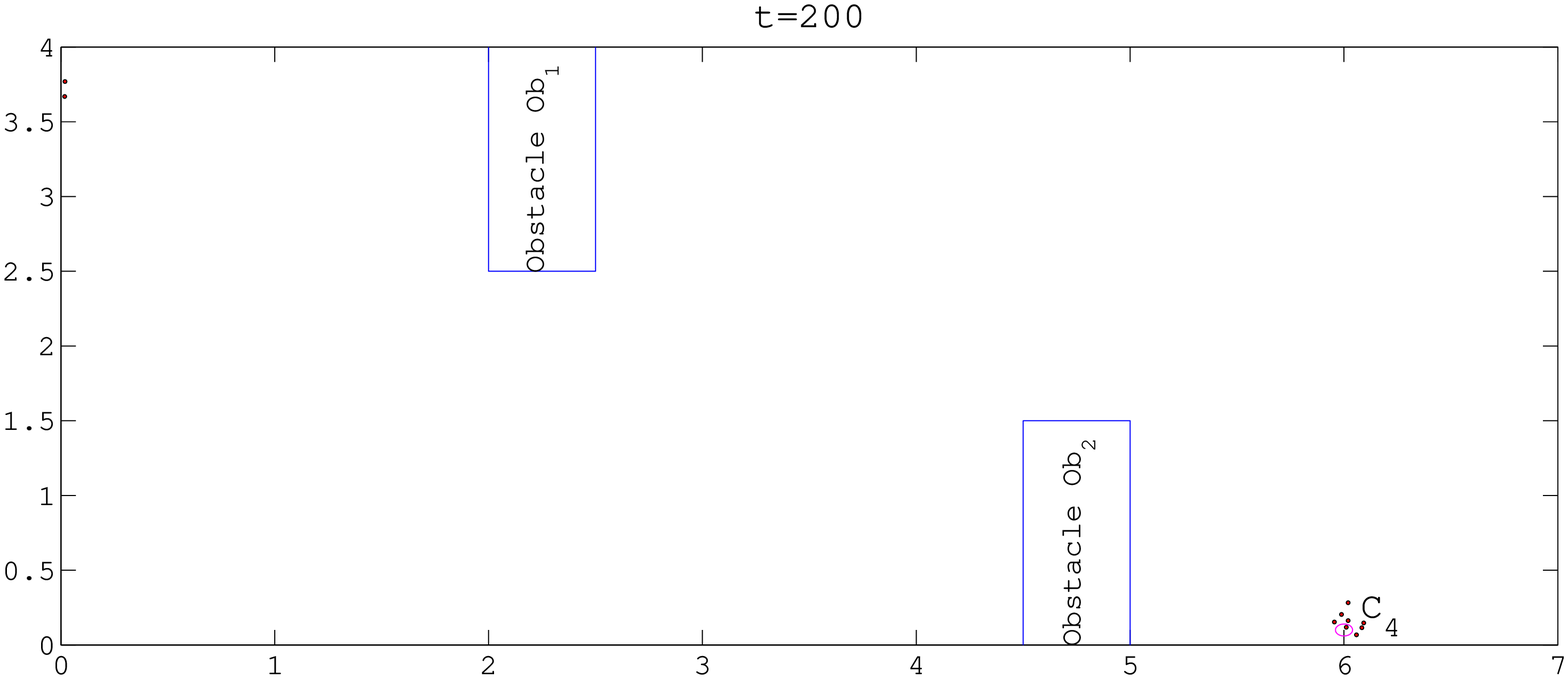} 
 \caption{A pattern of collective foraging. Positions of all ten fish in school are plotted at four  instants. The school reaches to the food resource at $C_4=(6, 0.1)$ while maintaining its school structure.} 
  \label{Config2_position}
 \end{center}
\end{figure}

In the last figure, Figure  \ref{Config2_position}, we give  a pattern of collective foraging. Positions of all ten fish in school are plotted at four instants $t=0, 50, 150, 200$. 


\begin{remark}
In numerical computations, we use the Euler explicit scheme for SDEs which has been introduced by Kloeden-Platen (\cite{Kloeden2005}). 
In general, this simple method has a strong order 0.5. Since the coefficients in our model are constant, the order is eventually one. Thus, this method is a suitable effective choice for our study as there are plenty of numerical simulations need to be made. 
\end{remark}

\section{Conclusions}  \label{section5}

We studied the  process of fish schools foraging in noisy environment with obstacle  from a mathematical point of view. We introduced the local rule {\rm (e)} for foraging and wrote out it into a mathematical formula. Then, we newly presented  the SDE model \eqref{eq4} describing the process by integrating the formula into our previous models. Our model described  the behavioral rules  of individual   precisely, and was tractable for mathematical treatments, especially performing  numerical simulations. 

Our numerical results qualitatively agreed with the interesting experimental observations  (\cite{Couzin1,Gotmark}) and empirical results  (\cite{Pitcher})  that the bigger the school size the larger the probability of foraging success. Our model, however, gave a prediction that this fact is not retained unboundedly. It is shown that there is an optimal value for school size  at which the probability of foraging success is highest. We may estimate it by means of numerical computations based on the model. Furthermore, the existence of this optimal value may be explained by the cohesiveness of school which is defined in \cite{LinhTonYagi}.

\section*{Acknowledgments} The authors heartily express their gratitude to the two anonymous reviewers for suggestions that greatly improved this manuscript.

\end{document}